\documentclass{article}
\usepackage[hypertex]{hyperref}
\usepackage{graphicx}
\usepackage{amsmath}
\usepackage{amssymb}
\usepackage{color}
\usepackage{mathrsfs}
\usepackage{verbatim}

\newcommand{\bbibitem}{\bibitem}

\newcommand{\llabel}[1]{{\label{#1}}}

\renewcommand{\r}[1]{(\ref{#1})}
\newcommand{\bi}{\begin{itemize}}
\newcommand{\ei}{\end{itemize}}
\newcommand{\bd}{\begin{description}}
\newcommand{\ed}{\end{description}}
\newcommand{\be}{\begin{enumerate}}
\newcommand{\ee}{\end{enumerate}}

\renewcommand{\i}{\item}

\newcommand{\bqn}{\begin{eqnarray}}
\newcommand{\eqn}{\end{eqnarray}}
\newcommand{\eqnn}{\nonumber\end{eqnarray}}
\newcommand{\eqnl}[1]{\llabel{#1}\end{eqnarray}}

\newcommand{\nn}{\nonumber}

\newcommand{\ba}[1]{\begin{array}{#1}}
\newcommand{\ea}{\end{array}}

\newcommand{\R}{\mathbb{R}}
\newcommand{\C}{\mathbb{C}}

\newcommand{\fine}{\end{document}}

\def \trait (#1) (#2) (#3){\vrule width #1pt height #2pt depth #3pt}
\def \qed{\hfill
        \trait (0.1) (6) (0)
        \trait (6) (0.1) (0)
        \kern-6pt
        \trait (6) (6) (-5.9)
        \trait (0.1) (6) (0)
\medskip}

\newtheorem{ml}{\bf Lemma}
\newtheorem{Theorem}{\bf Theorem}

\newtheorem{mrem}{\bf \underline{{\sl Remark}}}
\newtheorem{mcc}{\bf Corollary}
\newtheorem{Definition}{\bf Definition}
\newtheorem{mpr}{\bf Proposition}
\newtheorem{mproperty}{\bf Property}

\newcommand{\bt}{\begin{Theorem}}
\newcommand{\et}{\end{Theorem}}
\newcommand{\bl}{\begin{ml}}
\newcommand{\el}{\end{ml}}
\newcommand{\bp}{\begin{mpr}}
\newcommand{\ep}{\end{mpr}}
\newcommand{\bc}{\begin{mcc}}

\newcommand{\bproperty}{\begin{mproperty}}
\newcommand{\eproperty}{\end{mproperty}}

\newcommand{\ec}{\end{mcc}}
\newcommand{\bdeff}{\begin{Definition}}
\newcommand{\edeff}{\end{Definition}}
\newcommand{\brem}{\begin{mrem}\rm}
\newcommand{\erem}{\end{mrem}}

\newcommand{\proof}{{\bf Proof. }}
\newcommand{\ppotR}[3]
{

\begin{figure}\begin{center}
~\includegraphics[width=#3truecm]{./#1.eps}\\
\caption{#2}
\llabel{#1}
\end{center}
\end{figure}
\noindent$\!\!$}

\newcommand{\la}{\lambda}
\newcommand{\g}{\gamma}
\newcommand{\al}{\alpha}
\newcommand{\eps}{\varepsilon}
\newcommand{\de}{\delta}

\newcommand{\con}{{\cal C}}

\newcommand{\F}{{\cal F}}

\renewcommand{\O}{{\cal O}}

\renewcommand{\H}{{\cal H}}

\newcommand{\D}{{\cal D}}

\newcommand{\neigh}{neighborhood }

\newcommand{\eproof}{\hfill $\blacksquare$}

\newcommand{\und}{\underline}

\newcommand{\vep}{\varepsilon}

\newcommand{\ca}{c_{\alpha}}
\newcommand{\cac}{c^2_{\alpha}}
\newcommand{\sa}{s_{\alpha}}
\newcommand{\cs}{c_{s}}
\newcommand{\sis}{s_{s}}
\newcommand{\sac}{s^2_{\alpha}}
\newcommand{\sisc}{s^2_{s}}
\newcommand{\ra}{r(\alpha)}

\newcommand{\beq}{\begin{equation}}
\newcommand{\eeq}{\end{equation}}

\newcommand{\bpm}{\begin{pmatrix}}
\newcommand{\epm}{\end{pmatrix}}

\newcommand{\ga}{\gamma}

\newcommand{\dert}[1]{\frac{d}{d#1}}
\newcommand{\derp}[1]{\frac{\partial}{\partial #1}}

\newcommand{\NMON}{{\cal N}_{\mbox{\footnotesize mon}}}

\begin{document}

\title{{\LARGE{\sl {\bf Limit Time Optimal Synthesis\\ for a Control-Affine System on 
$S^2$}}}}
\author{P. Mason\thanks{Institut Elie Cartan UMR 7502, Nancy-Universit\'e/CNRS/INRIA POB 239,
54506 Vandoeuvre-l\`es-Nancy, France  {\tt Paolo.Mason@iecn.u-nancy.fr}},
R. Salmoni\thanks{Laboratoire des signaux et syst\`emes, Universit\'e Paris-Sud, CNRS, 
Sup\'elec, 91192 Gif-Sur-Yvette, France {\tt rebecca.salmoni@lss.supelec.fr}}, 
U. Boscain\thanks{SISSA, via Beirut 2-4 34014 Trieste, Italy {\tt boscain@sissa.it} and Le2i,  CNRS  UMR 5158, Universit\'e de Bourgogne, 9, avenue Alain Savary- BP 47870, 21078 Dijon, France}, 
Y. Chitour\thanks{Laboratoire des signaux et syst\`emes, Universit\'e Paris-Sud, CNRS, 
Sup\'elec, 91192 Gif-Sur-Yvette, France {\tt yacine.chitour@lss.supelec.fr}
\newline
The first author was (partially) supported by IdF--Aide au partage des projets
europ\'eens}
}

\maketitle

\vspace{.5cm} \noindent \rm

\begin{quotation}
\noindent  {\bf Abstract}
For $\al\in]0,\pi/2[$, let $(\Sigma)_\al$ be the control system $\dot{x}=(F+uG)x$,
where $x$ belongs to the two-dimensional unit sphere $S^2$, $u\in [-1,1]$
and $F,G$ are $3\times3$ skew-symmetric matrices generating
rotations with perpendicular
axes of respective length $\cos(\al)$ and $\sin(\al)$.
In this paper, we study the time optimal
synthesis (TOS) from the north pole $(0,0,1)^T$ associated to $(\Sigma)_\al$,
as the parameter $\al$ tends to zero. We first prove that the TOS
is characterized by a ``two-snakes'' configuration on the whole $S^2$,
except for a neighborhood $U_\al$ of the south pole $(0,0,-1)^T$ of diameter
at most $\O(\al)$.
We next show that, inside $U_\al$,
the TOS depends on the relationship between $r(\al):=\pi/2\al-[\pi/2\al]$ and
$\al$.
More precisely, we characterize three main
relationships, by considering sequences $(\al_k)_{k\geq 0}$ satisfying  $(a)$
$r(\al_k)=\bar{r}$; $(b)$ $r(\al_k)=C\al_k$ and $(c)$ $r(\al_k)=0$, where $\bar{r}\in (0,1)$ and $C>0$. 
In each case, we describe the TOS and provide, after a suitable rescaling,
the limiting behavior, as $\al$ tends to zero, of the corresponding TOS inside $U_\al$.
\end{quotation}




\newpage
\section{Introduction}
\llabel{intro}
Let $\al\in]0,\pi/2[$. On the unit sphere $S^2\subset\R^3$,
consider the control system $(\Sigma)_\al$ defined by
\beq
\label{cs-p}
(\Sigma)_\al\quad \dot{x}=(F+uG)x,~~x=(x_1,x_2,x_3)^T,~~~
\|x\|^2=1,~~|u|\leq1,
\eeq
where $F$ and $G$ are two  $3\times3$ skew-symmetric matrices
representing two orthogonal rotations with axes of length respectively
$\cos(\al)$ and $\sin(\al)$, $\al\in]0,\pi/2[$
(for the precise meaning of length, see Section~\ref{not0}).
With no loss of generality, we assume that
\beq
\label{FG0}
F:=\bpm
0&-\cos(\al)&0\\\cos(\al)&0&0\\0&0&0\epm\quad G:=\bpm 0&0&0\\
0&0&-\sin(\al)\\ 0&\sin(\al)&0\epm.
\eeq
In this paper, we aim at describing the time optimal
synthesis (TOS for short) from the north pole $N:=(0,0,1)^T$ for $(\Sigma)_\al$,
i.e. for every $\bar x\in S^2$ we want to find the time optimal
trajectory steering $N$ to $\bar x$ in minimum time (see Figure
\ref{f-inizio}).

In particular {\it we are interested in the qualitative shape of the
time optimal synthesis in a neighborhood of the south pole $S=(0,0,-1)^T$,
in the limit $\al\to0$.} The interest for that problem stems from motion planning
issues in aeronautics and quantum control, see \cite{y2,q5} for instance.
\ppotR{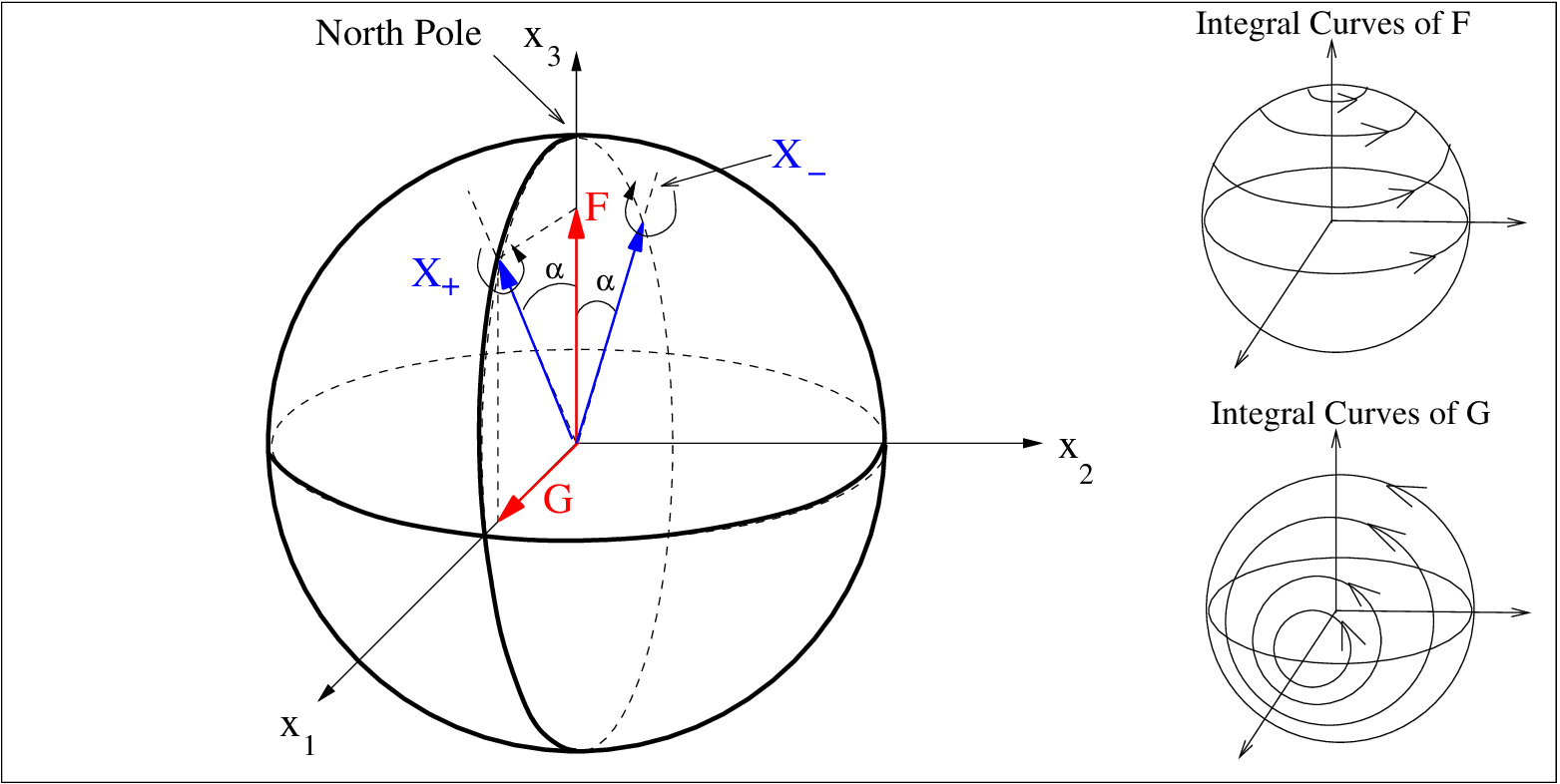}{
Geometric interpretation of the system
$(\Sigma)_\al$. The vector fields $X_+:=F+G$ and $X_-:=F-G$ are two
rotations of norm one making an angle $\al$ with the axis $x_3$.
}{12}


The present paper is actually a continuation of \cite{y2} in the sense that it answers
questions raised in the latter paper. There, the purpose was to
provide a lower and an upper bound for $N(\al)$, the maximum number of switchings for time
optimal trajectories for the left invariant control system
\bqn
\label{sysSO3}
(S)_\al\quad \dot g= g(F +uG), \quad g\in SO(3),\ \ |u|\leq 1,
\eqn
where $F$ and $G$ are defined in \r{FG0}.
Recall that, for such control systems, it is known (cf. for instance \cite{y2,q5})
that every time optimal trajectory is a finite concatenation of bang arcs
(i.e.  $u\equiv \pm 1$) or singular arcs ($u=0$). A bang arc is an integral trajectory
corresponding to the rotations
\bqn
X_+:=F+G,~~~X_-:=F-G,
\eqnl{X+-}
and is denoted by $e^{tX_\eps}x$,  $t\in[0,T]$, where $\eps=\pm$, $x$ is the starting point of
the bang arc and $T$ is its time duration.
Moreover, a switching time -- or simply a switching -- along a time optimal trajectory
is a time $t_0$ so that the control $u$ is not constant in any open neighborhood of $t_0$.

To estimate $N(\al)$, a suitable Hopf map $\Pi:SO(3)\rightarrow
S^2$ was introduced to
project $(S)_\al$ onto $(\Sigma)_\al$. In particular, every time optimal trajectory
of $(\Sigma)_\al$ is the projection by $\Pi$ of a time optimal trajectory
of $(S)_\al$. It results that, if a time optimal trajectory on $S^2$ has
a certain number of switchings, then this number is lower than or equal to
the maximum number of switchings for
the optimal problem on $SO(3)$.
The construction of time optimal trajectories of $(\Sigma)_\al$ was performed
according to the general theory of time optimal synthesis on 2-D manifolds developed in
\cite{ex-syn,automaton,tre,quattro,uno,due,sus1,sus2}, and recently gathered in the book
\cite{libro}.\\

The question of studying $N(\al)$ was first addressed in \cite{agra-sympl-x} where,
using the index theory developed by Agrachev, the authors proved that
$N(\alpha)\leq\left[\pi/\alpha\right]$, where
$\left[\cdot\right]$ stands for the integer part.  That result was not
only an indirect indication that $N(\alpha)$ would tend to infinity as
$\alpha$ tends to zero, but it also provided a hint on
the asymptotic of
$N(\alpha)$ as $\alpha$ tends to zero. Notice that for $\al=0$ the systems
\r{cs-p} and  \r{sysSO3} are not controllable.
With the techniques developed in \cite{y2}, enough properties for the TOS associated to
$(\Sigma)_\al$, $\al<\pi/4$, were identified in order to improve the
upper bound of \cite{agra-sympl-x} and to actually show that, for $\al$ small
$$
N(\alpha)\leq k_M+5,\quad\mbox{where}\quad k_M:=\left[\frac{\pi}{2\alpha}\right].
$$
In \cite{y2}, it is proved that, for $\al<\pi/4$, the
extremals associated to $(\Sigma)_\al$ (i.e. the trajectories candidate for time
optimality obtained after using the Pontryagin Maximum  Principle --PMP for short--), starting 
from the north pole $N$  are 
bang-bang trajectories,
i.e. finite concatenations of bang arcs of the type
$$e^{s_fX_{-{\eps'}}}e^{v(s_i)X_{\eps'}}\ldots e^{v(s_i)X_{-\eps}}e^{s_iX_\eps}N$$,
where the initial time duration $s_i$ verifies $s_i\in (0,\pi]$, all the time durations of the
interior
bang arcs are equal to $v(s_i)$, where the function $v$ is defined in
Eq.~\r{v()} below, and the final
time duration $s_f$ verifies $s_f\leq v(s_i)$. Of particular importance for the construction
of the TOS, are the \und{switching
curves}, i.e. the curves made by points where the control
switches from $+1$ to $-1$ or viceversa and defined inductively by
\beq\label{swi0} C_1^{\eps}(s)=
e^{X_\eps v(s)} e^{X_{-\eps}s}  N,\ \ \
C^{\eps}_{k}(s)=e^{X_\eps v(s)}C^{-\eps}_{k-1}(s),  \mbox{ (where
$\eps=\pm1$ and $k=2,....,k_M$).} 
\eeq
Since the PMP gives just a necessary condition for optimality, 
it is crucial to determine the time after which an extremal is no more optimal. 
In \cite{y2}, we showed that the number of bangs 
must be
lower than or equal to $k_M+1$ and the extremals cover the sphere $S^2$ according to the
''two-snakes''
configuration as depicted in Figure \ref{f-serpentone2}.  The two ``snakes'' correspond to 
extremal trajectories starting respectively with control $+1$ and $-1$. For more details,
see \cite{y2}.
\ppotR{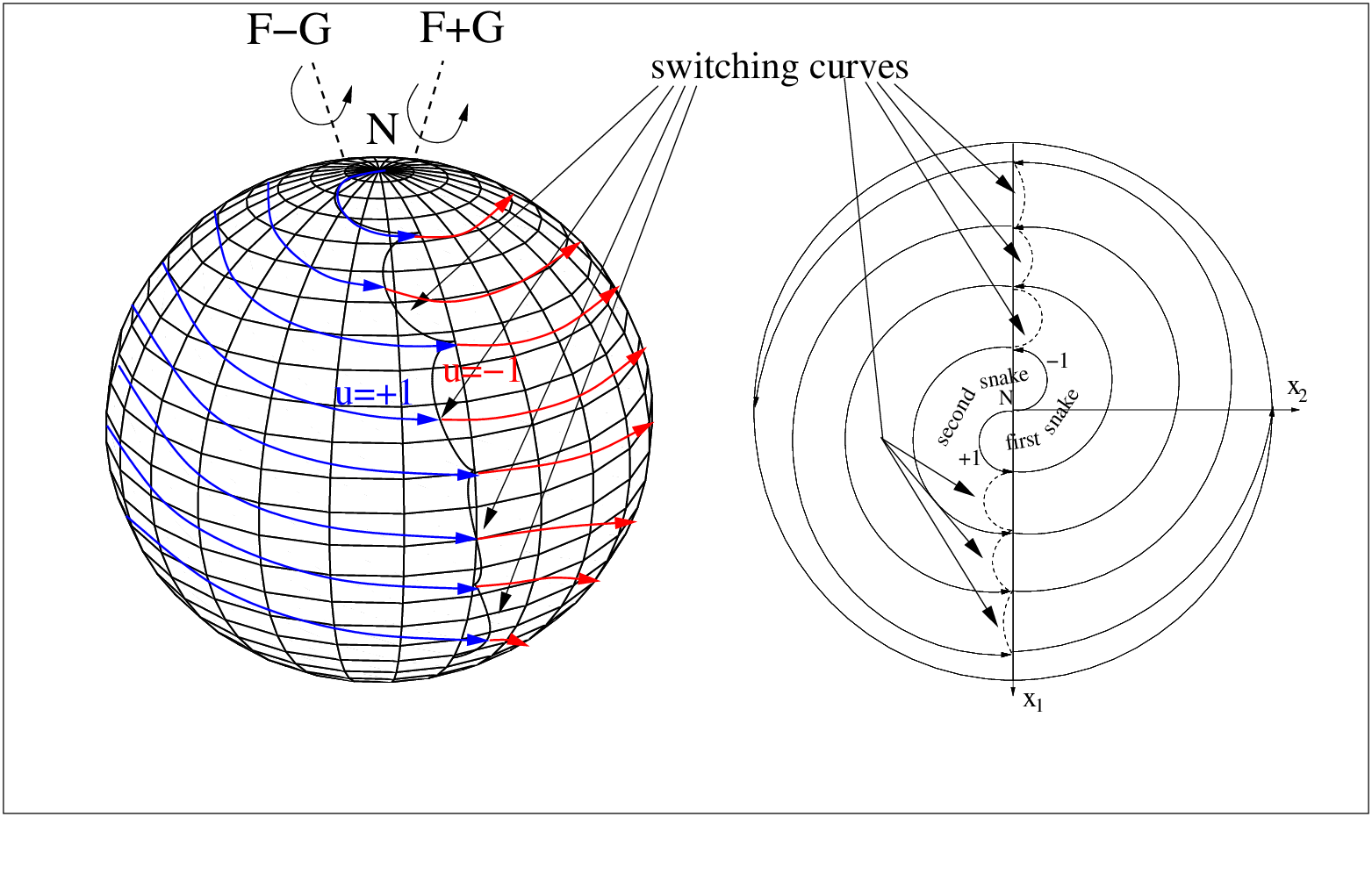}{The ``two-snakes'' configuration defined by the extremal flow. Notice 
that this set of trajectories covers the whole sphere, but in principle not all extremals 
are optimal and a point can be reached by more than one trajectory at the same or at 
different times.}{13}

\begin{figure}
\begin{center}
\includegraphics[width=15truecm,angle=0]{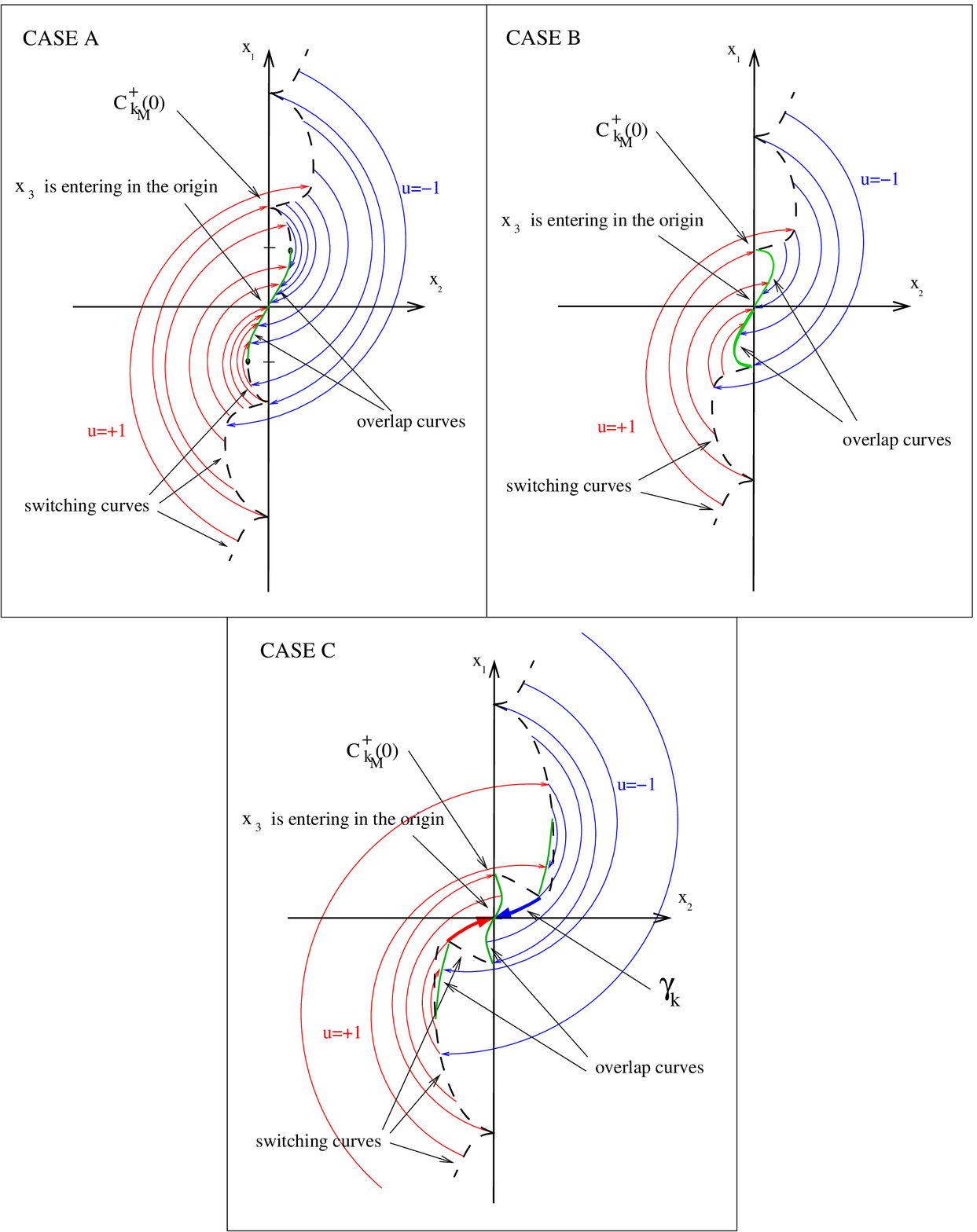}
\caption{Conjectured shapes of the synthesis in a neighborhood of the south pole.
Switching curves are $\con^1$ curves made by points in which the control 
switches from $+1$ to $-1$ or viceversa. Overlap curves are $\con^1$ 
curves made by points reached optimally by more than one trajectory.
The curve $\g_k$ is a bang arc that is also an overlap curve since 
trajectories having a different history travel on it at the same time.
} \label{f-TUTTALASINTESI}
\end{center}
\end{figure}

However, in \cite{y2}, we were not able to
construct the complete TOS associated to $(\Sigma)_\al$. In particular, we could not show the
optimality of all the extremals up to $k_M-1$ bangs arcs and we could not complete
analytically the construction of the synthesis in a
neighborhood of the south pole $S$. There, the minimum time front develops singularities due
to the compactness of $S^2$. We only provided numerical simulations describing the evolution of the
extremal front in a neighborhood of the south pole.
As $\al\to0$, these numerical simulations suggested the emergence of an
interesting phenomenon (see Fig.~\ref{f-TUTTALASINTESI}) : define the remainder
\bqn
r(\al):=\pi/2\al-[\pi/2\al].
\eqnl{remainder}

Then, there are three possible patterns of TOS in the neighborhood of the south pole $S$,
each of them depending on a relation between $r(\al)$ and $\al$. See 
Section~\ref{des000} below, where these relations are
formulated as well as conjectures.\\

In \cite{q5}, the TOS for $(\Sigma_\al)$ was studied in the context of 
quantum control. The
control system $(\Sigma_\al)$ describes the population transfer problem for the $x_3$-component of
the spin of a (spin
1/2) particle, driven by a magnetic field, which is constant along the $x_3$-axis and controlled
along the $x_1$-axis, with bounded amplitude. In that paper, the TOS, for $\al\geq\pi/4$ was completed
and, in the case  $\al<\pi/4$, further information was obtained, for what concerns time optimal
trajectories steering the north to the south pole
(in fact the most interesting trajectories for the quantum mechanical problem).
Such optimal trajectories belong to a set  $\Xi$ containing at most  8 trajectories, half of them
starting with control $+1$ and the other half starting with control $-1$, and switching exactly at
the same times. It was also proved that the cardinality of $\Xi$ depends on the remainder $r(\al)$
defined in Eq.~\r{remainder}. For instance, for $\al$ and  $r(\al)$ small enough, then $\Xi$
contains exactly $8$ trajectories (four of them are optimal)
while if  $r(\al)$ is
close to $1$, then $\Xi$ contains only $4$ trajectories (two of them are optimal).

The purpose of the present paper consists in studying the TOS
associated to $(\Sigma)_\al$ as $\al$ tends to zero, focusing in particular on its
behavior inside a neighborhood of the south pole. 
Roughly speaking, we want to determine, as $\al$ tends to zero, what could be a possible
limit for the TOS associated to $(\Sigma)_\al$ (as suggested for instance by the patterns depicted in
Fig.~\ref{f-TUTTALASINTESI}) and then to prove the convergence (in some suitable sense)
of the TOS associated to $(\Sigma)_\al$ to that limit. To proceed, we embark on the study of
a geometric object $\F(\al,T)$ called the {\it extremal front at time $T$} along $(\Sigma)_\al$
and defined as the set of points reached at time $T$ by extremal
trajectories starting from $N$ (see Section~\ref{EXTRE} for a precise definition). The extremal 
front $\F(\al,T)$ contains 
the {\it minimum time front} $OF(\al,T)$, i.e.
the set of points reached at time $T$ by time optimal trajectories. 
When $\F(\al,T)=OF(\al,T)$,
we say that $\F(\al,T)$ is {\it optimal}.


We first prove, in the case in which $k_M$ is odd (being the other case analogous),
that the extremal front  $\F(\al,k_M\pi)$ is
made up of the union of two curves
$\mathcal{E}^\eps(\al,\cdot):(0,\pi]\rightarrow S^2$, $\eps=\pm$, 
with $\mathcal{E}^\eps(\al,\cdot)=\Pi_{x_3}\mathcal{E}^{-\eps}(\al,\cdot)$, 
where $\Pi_{x_3}$ is the orthogonal symmetry with respect to the $x_3$-axis. 
Moreover, for $\al$ small enough, $\mathcal{E}^\eps(\al,\cdot)$ admits a 
convergent power series of the 
type $\sum_{l\geq 0}f_l^\eps(s,r(\al))\al^l$, where
the $f_l^\eps(s,r)$ are real-analytic functions of $(s,r)\in\R^2$, 
$2\pi$-periodic in $s$ with
\beq\label{psw2}
f_0^+(s,r)=\bpm 0\\ 0\\-1\epm, \qquad f_1^+(s,r)=\bpm -2r\cs\\
2r\sis\\0\epm, \qquad f_2^+(s,r)=\bpm
\frac{\pi}2(4r+\cs)\sisc\\\frac{\pi}4(3+8r\cs+c_{2s})\sis\\2r^2\epm.
\eeq
As a trivial consequence, we deduce that for $r\in [0,1]$, $s\in \R$ and $\al$ small enough,
we have 
\bqn
\mathcal{E}^\eps(\al,s)=f_0^\eps(s,r(\al))+f_1^\eps(s,r(\al))\al+f_2^\eps(s,r(\al))\al^2+\O(\al^3).
\llabel{yacineseiunbalordo}
\eqn
and 
\bqn
\derp{s}\mathcal{E}^\eps(\al,s)=\derp{s}f_1^\eps(s,r(\al))\al+
\derp{s}f_2^\eps(s,r(\al))\al^2+\O(\al^3).
\llabel{derfronte}
\eqn
where $|\O(\al^3)|\leq \bar C |\al|^3$ with $\bar C>0$ constant independent of $(r,s,\al)$. 


Then we show that  $\F(\al,T)$ is actually optimal
for $T\leq (k_M-1)\pi$ and $\alpha$ small
enough (see Remark~\ref{opt0} below). 
Moreover, we show that $\F(\al,(k_M-1)\pi)$
is a circle of radius $2(1+\ra)\alpha$ up to order $\al^2$
(see Remark~\ref{lem-est} below). 
As a consequence of the optimality of  $\F(\al, (k_M-1)\pi)$, we get that 
all the extremals of the ``two-snakes'' configuration depicted in Figure
\ref{f-serpentone2} are optimal up to time $(k_M-1)\pi)$. In other words, if $U_\al$ 
is the connected 
component of $S^2\setminus \F(\al,(k_M-1)\pi)$ containing the south pole, we  obtain the
 optimal synthesis on $S^2\setminus U_\al$. Notice that  $U_\al$ is a  neighborhood  of the south
pole of size proportional to $\alpha$.



The expressions \r{yacineseiunbalordo}--\r{derfronte} are  central tools to understand the 
possible asymptotic behaviors of the TOS associated to $(\Sigma)_\al$, as $\al$ tends to zero. 

For this purpose we observe that the expressions of  $f^+_1$ and $f^+_2$ in \r{psw2}, 
depend explicitly on the remainder $r(\al)$. This fact suggests the need to impose 
particular relationships between $\al$ and $r(\al)$ in order to define any
asymptotic behavior. 
In other words we must let $\al$ goes to zero, only
along certain subsequences $(\al_k)_{k\geq 0}$ where a specific
relationship holds between $\al_k$ and $r(\al_k)$. The analysis of
Eq.~\r{psw2} will help us to determine such relationships and to prove 
that the conjectures made in \cite{y2} about the qualitative shape of the synthesis 
near the south pole were true (see Section~\ref{des000} and Figure~\ref{f-TUTTALASINTESI}). 
In particular we will see that there are exactly three qualitatively different 
asymptotic behaviours of the synthesis as $\al$ goes to zero, described 
by the following cases.

First, we analyze the case in which $\al$ is arbitrarily small, with $r(\al)\in(0,1)$ 
uniformly far from $0$ and $1$.
To simplify further the discussion, it is reasonable to consider the following.
\begin{center}
(C1)~~~For $\bar{r}\in (0,1)$, let $\al$ tend to zero along the subsequence 
 $\al_k:=\frac{\pi}{2 (k+\bar r)}$, so that $r(\al_k)=\bar{r}$.
\end{center}
In this case $\mathcal{E}^\eps (\al,\cdot)$ is approximated, up to order $\al^2$, by the expression
$S+f_1^\eps(\cdot,\bar r)$. As a consequence 
$\F(\al,k_M\pi)$ is 
approximately a circle of radius $2\bar r\al$ centered at the south pole.
We are then able to give a qualitative description of the optimal synthesis, as stated below 
in Theorem~\ref{the1}. This synthesis  turns out to be exactly the one described in Figure 
~\ref{f-TUTTALASINTESI} (case B), as  predicted in \cite{y2}.

It remains then to study the cases in which $r(\al)$ can be arbitrarily close to $0$ or $1$.
For this purpose we first consider the case in which $r(\al)/\al$ remains
bounded above and below by positive constants as $\al$ tends to zero.
From Eq.~\r{psw2} it is clear that this is equivalent to say that $f_1^\eps(\cdot,r)\al$
is comparable to $f_2^\eps(\cdot,r)\al^2$.
For simplicity we consider the following.
\begin{center}
(C2)~~~For $C>0$, let $\al$ tend to zero along a subsequence $(\al_k)_{k\geq 0}$ 
such that $r(\al_k)=C\al_k$.
\end{center}
In this case
$\mathcal{E}^\eps(\al,\cdot)$ is well approximated by $S+(f_1^\eps(\cdot,C)+
f_2^\eps(\cdot,0))\al^2$.
If $C>\pi/4$, the  synthesis is equivalent the one of the  previous case. On the other hand 
if  $C<\pi/4$  the  synthesis is more complicated (see Section \ref{Kal}) and it   
turns out to be exactly the one described in Figure~\ref{f-TUTTALASINTESI} (case C), as  
predicted in 
\cite{y2}.

If $\al$ and $r(\al)$ tend to zero with $r(\al)/\al$ tending  to infinity (resp. to zero) it 
is possible to 
see that  the synthesis is qualitatively equivalent to the one of case (C1) (resp. (C2)). 

The third interesting case is the following.
\begin{center}
(C3)~~~Let $\al$ tend to zero along the subsequence $\al_k:=\frac\pi{2 k}$, so that $r(\al_k)=0$.
\end{center}
In this case the extremal front at time $k_M\pi$ contains the south pole and 
the corresponding optimal front reduces to that point. The optimal synthesis
is then described starting from the extremal front $\F(\al,(k_M-1)\pi)=OF(\al,(k_M-1)\pi)$, and it 
corresponds to the one  described in Figure~\ref{f-TUTTALASINTESI} (case A), as  predicted in
\cite{y2}.

With similar arguments, one can see that in the case in which $\al$ is small and $r(\al)$ 
is close to $1$,
the optimal synthesis is qualitatively equivalent either to that of Case (C1) or to that of 
Case (C3), and this concludes the description of the possible asymptotic behaviors as $\al$
tends to $0$.

\brem
It is interesting to notice that numerical 
simulations show that for $\al$ decreasing to 
zero continuously, the qualitative shape of 
the optimal synthesis described in 
Figure~\ref{f-TUTTALASINTESI} alternates cyclically in the order BCABCA....
 \erem

Let us describe the results obtained in the case (C1) in more details.
Since $\F(\al,k_M\pi)$ is approximated,
up to $\O(\al^2)$, by a circle of center $S$ and radius
$2\bar{r}\al$, we are able to show that it is optimal, so that all the extremals of the 
``two-snakes'' 
configuration depicted in Figure~\ref{f-serpentone2} are optimal up to time $k_M\pi$. 
In other words, if $V_\al$ is the connected component of $S^2\setminus \F(\al,k_M\pi)$ 
containing the south pole, we  obtain the optimal synthesis on $S^2\setminus V_\al$. 

As $\al$ tends to zero, $V_\al$ collapses on $S$. Hence one must rescale
the problem by a factor $1/\al$, in order to describe the TOS
inside $V_\al$. 
Also notice that since we are in a neighborhood of the south pole we can project the problem on the plane $(x_1,x_2)$.
We are now in a position to define a possible limit behavior for the TOS
inside $V_\al$. Let $M_{\al}$ be the linear mapping from $\R^3$
onto $\R^2$ defined as the composition of the projection
$(x_1,x_2,x_3)\mapsto (x_1,x_2)$ followed by the dilation by
$1/\al$. Denote by $(\widetilde{\Sigma})_{\al}$ (resp. $\widetilde{OF}(\al,k_M\pi)$) the image by $M_{\al}$
of $(\Sigma)_\al$ (resp. $OF(\al,k_M\pi)$). Then,
$(\widetilde{\Sigma})_{\al}$ is a 
perturbation by $\O(\al^2)$ of the standard linearized
pendulum 
\beq\label{lineare0} (Pen):\quad\left\{ \ba{l}
  \dot{z_1}=-  z_2, \\
  \dot{z_2}= z_1 +u,\ \ \ \ (z_1,z_2)\in\R^2,\ \ |u|\leq 1\ea \right. \\
\eeq
while $\widetilde{OF}(\al,k_M\pi)$ is a perturbation by $\O(\al^2)$ of 
$C(0,2\bar{r})$, the planar circle of center $(0,0)$ and
radius $2\bar{r}$. As a consequence, the candidate limit TOS inside
$V_\al$ is the one associated to the problem of reaching  in minimum 
time every point of  the ball $B(0,2\bar{r})$ starting  
from $C(0,2\bar{r})$, along the
dynamics of the standard linearized pendulum. To prove such a
result, we first study the above mentioned optimal control problem
and show that the corresponding TOS is
characterized by an overlap curve $\gamma^o_{pen}$, which is the set 
of points $z\in\R^2$ with $z_1z_2\geq 0$ and belonging to the locus (see Figure \ref{f-over-pen})

$$z_1^4+z_2^4+2z_1^2z_2^2-4\bar r^2z_1^2+(4-4\bar r^2)z_2^2=0\,.$$
The optimal synthesis inside $C(0,2\bar{r})$ is then described by the 
following feedback, defined on $B(0,2\bar{r})\setminus\gamma^o_{pen}$: 
``above'' $\gamma^o_{pen}$, the control $u$ is constantly equal to $-1$
and ``below'' $\gamma^o_{pen}$, it is constantly equal to $1$ (see
Fig.~\ref{sintesilineare}). 
Finally, the asymptotic result we prove
in Section~\ref{s-rwlp} is the following.

\ppotR{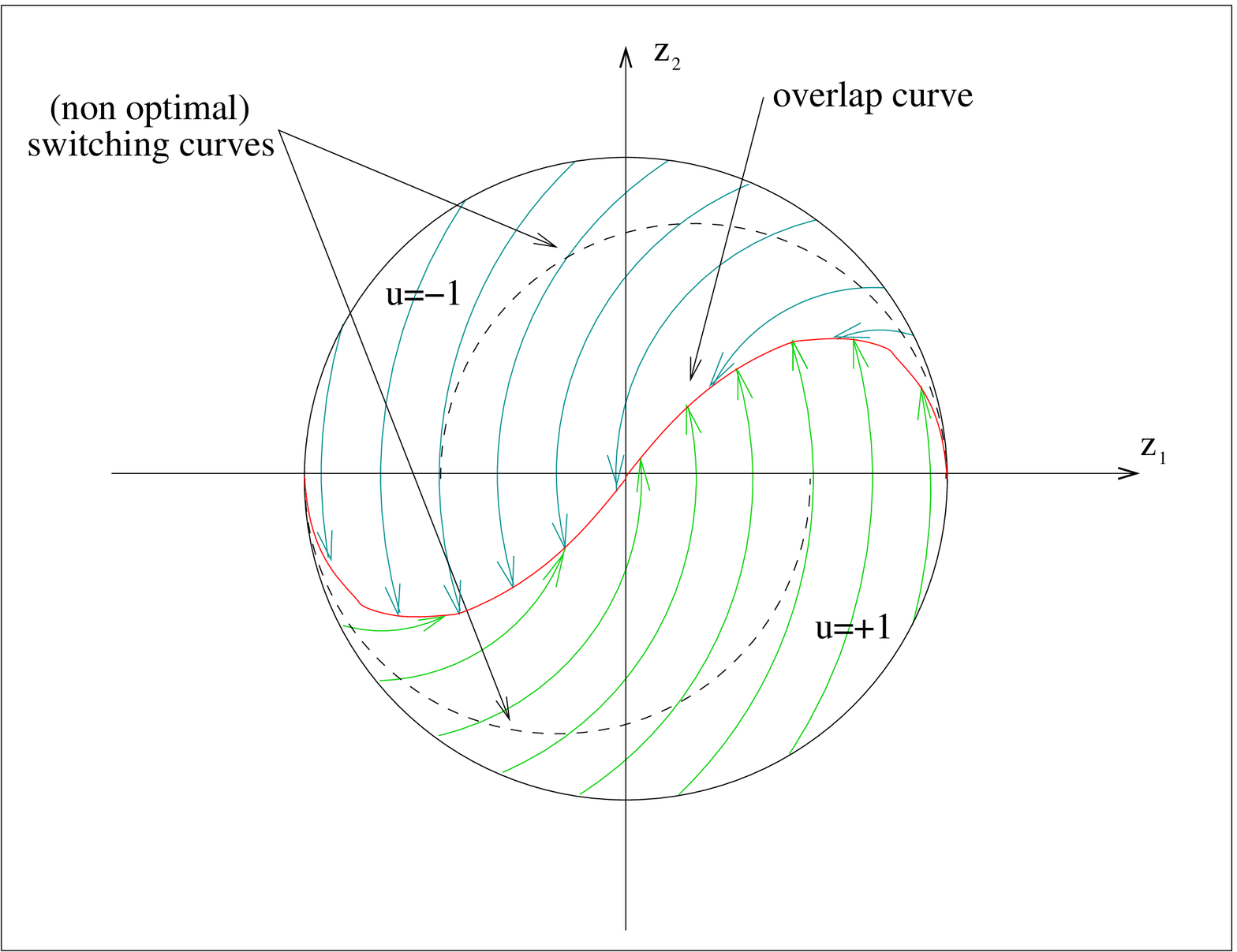}{Optimal synthesis for the linear pendulum}{10}

\bt\label{the1} 
For $\bar{r}\in (0,1)$, let $(\al_k)_{k\geq 1}$ be the sequence defined by
$\al_k:=\frac{\pi}{2(k+\bar r)}$ for $k\geq 1$. Consider $\gamma^o_{pen}$,
the overlap curve of the TOS for the
optimal control problem consisting of starting from $C(0,2\bar r)$, the
planar circle of center $(0,0)$ and radius $2\bar r$, and reaching in
minimum time every point of $B(0,2\bar r)$ along the control system
\r{lineare0}. 
Then, for $k$ large enough, the
TOS associated to $(\widetilde{\Sigma})_{\al_k}$ inside
$\widetilde{OF}(\al_k,k_M\pi)$ is characterized by an overlap curve
$\gamma^o_{\al_k}$ so that 
the optimal feedback takes the
value $-1$ ``above'' $\gamma^o_{\al_k}$, and the value $1$ ``below'' 
$\gamma^o_{\al_k}$. Moreover, $\gamma^o_{\al_k}$ converges to
$\gamma^o_{pen}$ in the $C^0$ topology, uniformly with respect to $\bar r$ in
any compact interval of $(0,1)$.
\et
The results in the cases (C2) and (C3) are described in more details in Sections 
\ref{Kal} and \ref{r=0}.

\brem
Notice that the sequence $(\al_k)_{k\geq 1}$ defined above has been
chosen in order to simplify the previous statement. Indeed the same result 
could be restated in a more general way by taking an arbitrary sequence
$(\tilde{\al}_k)_{k\geq 1}$ converging to zero and such that $r(\tilde\al_k)$
converges to $\bar r$, or letting the remainder vary on a compact subinterval
of $(0,1)$.
\erem


The paper is organized as follows. In the second section, we collect basic facts, notations, 
results, and conjectures of \cite{y2}. The third section
gathers the detailed description of the extremal front and the proof of Eq.~\r{psw2}.
Sections~\ref{barr}, ~\ref{Kal} and ~\ref{r=0} treat respectively the cases (C1), (C2) and (C3).
In the appendix, we finally prove 
a technical result needed throughout the paper.

\section{Notations and previous results}
\subsection{Basic Facts}
\label{ss-basic}

\bdeff
An admissible control $u(.)$ for the system \r{cs-p}--\r{FG0} is a
measurable function $u(.): [a,b]\to [-1,1]$,
while an admissible trajectory is an absolutely continuous function
$x(.):[a,b]\to S^2$ satisfying \r{cs-p} a.e.
for some  admissible control $u(.)$.
If $x(.)$ is an admissible trajectory and $u(.)$ the corresponding
control, we say that $(x(.),u(.))$ is an
admissible pair.
\edeff
For every $\bar x\in S^2$,  the minimization problem consists of determining an
admissible pair steering the north pole to $\bar x$ in minimum time.
More precisely\\\\
\medskip\noindent
{\bf Problem (P)}
{\it Consider the control system
\r{cs-p}-\r{FG0}.
For every $\bar x\in S^2$, find an admissible pair
$(x(.),u(.))$ defined on $[0,T]$ such that
$x(0)=N$, $x(T)=\bar x$ and $x(.)$ is time optimal.
}\\\\
An optimal synthesis from the north pole (in the following optimal synthesis, for short) 
is the collection of all the solutions to the problem {\bf (P)}.
More precisely
\bdeff {\bf (Optimal Synthesis)} An optimal synthesis for the problem
{\bf (P)} is the
collection of all time optimal trajectories
$\Gamma=\{x_{\bar x}(.):[0,b_{\bar x}]\mapsto S^2 $, $\bar x\in S^2
:~~x_{\bar x}(0)=N,
~x_{\bar x}(b_{\bar x})=\bar x\}$.
\edeff
For more elaborated definitions of optimal synthesis
see \cite{libro,piccoli-sussmann} and
references therein.
The standard tool to look for optimal trajectories is a first order
necessary condition for optimality  known as the  Pontryagin Maximum Principle (PMP for short),
(cf. \cite{agra-book,pontlibro})
as stated below for our minimum time problem on $S^2$.\\

Define the following real-valued map on $T^*S^2\times [-1,1]$, called
Hamiltonian,
$$
\H(\la,x,u):=<\la, (F+uG)x>.
$$
Set:
\bqn
H(\la,x):=\max_{v\in [-1,1]}\H(\la,x,v).
\eqnl{brrrr}
The PMP asserts that, if
$\gamma:[a,b]\to S^2$ is
a time optimal trajectory corresponding to a control $u:[a,b]\to [-1,1]$, then
there exists {\sl a nontrivial field of covectors along} $\ga$, that is
a never vanishing absolutely continuous function
$\la:t\in [a,b]\mapsto \la(t)\in T^\ast_{\ga(t)} S^2$  and a constant
$\la_0\leq 0$ such that, for a.e. $t\in Dom(\g)$,
we have:
\begin{description}
\item[i)]  $\dot \la(t)=-\frac{\partial \H}{\partial
x}(\la(t),\gamma(t),u(t))=-\la(t)(F+u(t)G)$,
\item[ii)]  $\H(\la(t),\ga(t),u(t))+\la_0=0$,
\item[iii)] $\H(\la(t),\ga(t),u(t))=H(\ga(t),\la(t))$.
\end{description}
In the more general case in which the target and the initial datum (also called \emph{source}) are two smooth
manifolds $\mathcal{N}_0$ and $\mathcal{N}_1$ the previous statement must be 
modified by adding the so-called \emph{transversality conditions}:
\bd
\i[iv)] $<\la(a),v>=0\quad\forall v\in T_{\gamma(a)}\mathcal{N}_0\,,\quad<\la(b),w>=0\quad\forall w\in T_{\gamma(b)}\mathcal{N}_1$.
\ed
\brem
\llabel{r-postPMP}
A trajectory $\g$
(resp.  a couple $(\g,\la)$)
satisfying the conditions given by the PMP is said to be an {\sl
extremal} (resp.  an {\sl extremal
pair}).
An extremal  corresponding to $\la_0=0$ is said to be an {\sl
abnormal
extremal}, otherwise we call it a {\sl normal extremal}.
\erem
\bdeff {\bf (bang, singular for the problem  \r{cs-p}-\r{FG0})}
\llabel{d-BS-1}
A control $u(.):[a,b]\to[-1,1]$  is said to be a \und{bang} control if
$u(t)=+1$
a.e. in  $[a,b]$ or  $u(t)=-1$ a.e.  in  $[a,b]$. A control
$u(.):[a,b]\to[-1,1]$ is said to be a
\und{singular} control if $u(t)=0$, a.e. in  $[a,b]$. A finite concatenation
of bang controls is
called a \und{bang-bang} control. A \und{switching time} of $u(.)$ is a
time
$\bar t\in[a,b]$ such that, for every $\vep>0$, $u$ is not bang or
singular on
$(\bar t-\vep,\bar t+\vep)\cap
[a,b]$.
A trajectory of the control system
\r{cs-p}, \r{FG0} is said to be a bang
trajectory (or arc), singular trajectory (or arc), bang-bang
trajectory, if it
corresponds respectively
to a bang control, singular control, bang-bang control. If $\bar t$ is a switching time, the
corresponding point on
the trajectory $x(\bar t)$ is called a \und{switching point}.
\edeff
\subsection{Description of Previous Results}\label{des000}
In  \cite{y2,q5} it was proved that,
for every couple of points, there exists a time
optimal trajectory joining them.
Moreover it was proved that
every time optimal trajectory is a finite
concatenation of bang and singular trajectories.
More precisely we have:
\bp
For the minimum time problem associated to 
\r{cs-p}-\r{FG0},
for each pair of points $p$ and $q$ belonging to $S^2$,
there exists a time optimal trajectory joining $p$ to $q$.
Moreover every  time optimal trajectory for \r{cs-p}-\r{FG0} is a
finite concatenation of bang and singular trajectories.
\ep
Notice that the previous proposition does not apply
if $\al=0$ or $\al=\pi/2$, since in these cases the controllability property is lost.

In \cite{q5} it has been proved that $\al=\pi/4$ is a bifurcation for the
qualitative shape of the time optimal synthesis, for instance the
time optimal synthesis contains a singular arc if and only if $\al>\pi/4$.
Since in this paper we are interested in the limit $\al\to0$, in the
following we always assume $\al<\pi/4$.
In this case, using the PMP, the following
properties characterizing the optimal trajectories were established in \cite{y2}:
\bd
\i[\bf i)] $x(.)$ is bang bang;
\i[{\bf ii)}] the
duration $s_i$ of the first bang
arc satisfies $s_i\in(0,\pi]$,
\i[{\bf iii)}] the time duration
between two consecutive switchings is the same for all \underline{interior
bang arcs} (i.e. excluding the first and the last bang)
and it is equal to $v(s_i)$, where $v(\cdot)$ is the following function 
\bqn
v(s)=\pi+2 \arctan\left(\frac{\sin(s)}
{\cos(s)+\cot^2(\alpha)}\right).
\eqnl{v()}
One can immediately check that this function satisfies $v(0)=v(\pi)=\pi$
and $v(s)>\pi$ for every $s\in(0,\pi)$,
\i[{\bf iv)}] the time duration of the last arc is $s_f\in(0,v(s_i)]$,
\ed
Moreover, thanks to the analysis given in \cite{y2}, one easily gets
(always in the case $\al<\pi/4$):
\bd
\i[{\bf v)}] the number of switchings $N_x$ of $x(.)$ satisfies the
following inequality
\bqn
N_x\leq\left[ \frac{\pi}{2\al}\right]+1.
\eqnl{eq-nmax}
\ed
\ppotR{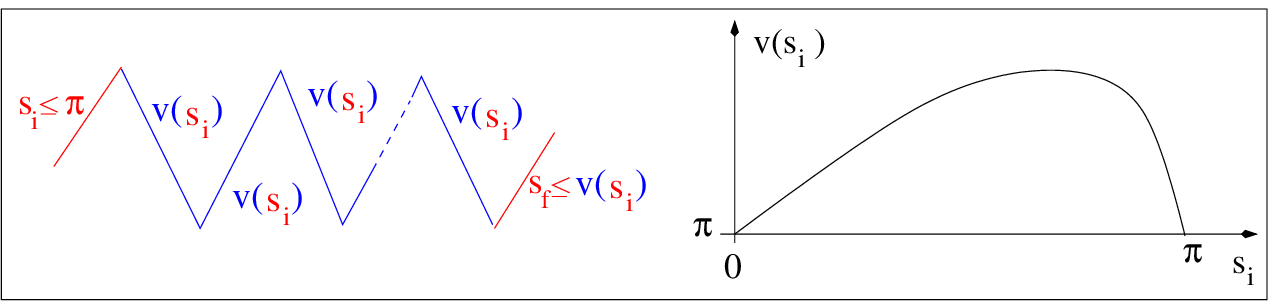}{Time optimal trajectories for $\al<\pi/4$}{12}
Conditions {\bf i)}-{\bf v)} define a set of
candidate optimal trajectories. Notice that conditions {\bf i)}-{\bf v)}
are just
necessary conditions for optimality and one is faced with the problem of
selecting, among them, those that are really
optimal. In particular, given a trajectory satisfying  conditions {\bf
i)}-{\bf v)}, one would like to find the time
after which it is no more optimal.




Some questions remained unsolved, in particular questions relative
to local optimality of the switching curves defined in Eq.~\r{swi0}.
Roughly speaking we say
that a switching curve is locally optimal if it never ``reflects''
the trajectories (see Fig.~\ref{f-nonlocop} A).\footnote{ More
precisely consider a smooth switching curve $C$ between two smooth
vector field $Y_1$ and $Y_2$ on a smooth two dimensional manifold.
Let $C(s)$ be a smooth parametrization of $C$. We say that $C$ is
\underline{locally optimal} if, for every $s\in Dom(C)$, we have
$\dot C(s)\neq\al_1 Y_1(C(s))+\al_2 Y_2(C(s)), \mbox{ for every
}\al_1,\al_2\mbox{ s.t. }\al_1\al_2\geq0.$ The points of a switching
curve on which this relation  is not satisfied are usually called
``conjugate points''. See Fig.~\ref{f-nonlocop}. }
When a family of trajectories is reflected by a switching curve then
local optimality is lost and some \und{cut locus} appear in the
optimal synthesis.

\bdeff A cut locus is a set of points reached at the same time by
two (or more) time optimal trajectories. A subset of a cut locus
that is a connected $\con^1$ manifold is called \und{overlap curve}.
\edeff

An example showing how a ``reflection'' on a switching curve
generates a cut locus is portrayed in Fig.~\ref{f-nonlocop} B and C.
More details are given later. More precisely the following questions
were formulated in \cite{y2}: \ppotR{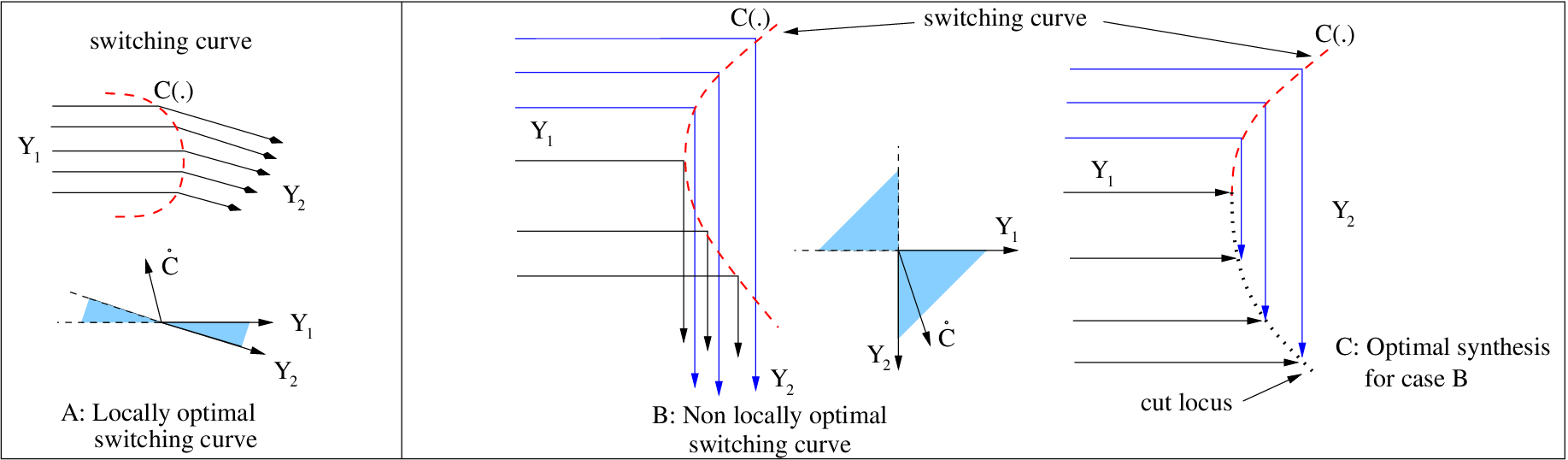}{Locally optimal switching curves
and non locally optimal switching curves with the corresponding
synthesis}{16.5}

\bd \i[Question 1] Are the switching curves $C^\eps_k(s)$,
$s\in(0,\pi]$, locally optimal? More precisely, one would like to
understand  how the candidate optimal trajectories described above
lose optimality.

\i[Question 2] What is the shape of the optimal synthesis in a
\neigh of the south pole? \ed
Numerical simulations suggested some conjectures regarding the above
questions. More precisely, in \cite{y2} the following conjectures were made:
\bd 
\i[C1] The curves
$C^\eps_k(s)$, ($k=1,...,k_M$) are locally optimal if and only if
$k\leq \left[\frac{\pi-\al}{2\al}\right]-1$.
\ed
Analyzing the evolution of the minimum time wave front in a \neigh
of the south-pole, it is reasonable to conjecture that:
\bd \i[C2] The shape of the optimal synthesis in a neighborhood of
the south pole depends on the remainder $r(\al)$ defined in Eq.~\r{remainder}.
Notice that $r(\al)$ belongs to the interval $[0,1)$. More precisely,
it was conjectured in \cite{y2} that for $\al\in(0,\pi/4)$, there exist two positive
numbers $\al_1$ and $\al_2$ such that $0<\al_1<\al<\al_2<2\al$ and:
\bd \i{}\underline{CASE A: $\ra\in (\frac{\al_2}{2\al},1)$}. The switching curve
$C^\eps_{k_M}$ glues to an overlap  curve that passes  through the
origin (Fig.~\ref{f-TUTTALASINTESI}, Case A).

\i{}\underline{CASE B:
  $\ra\in[\frac{\al_1}{2\al},\frac{\al_2}{2\al}]$}.
The switching curve $C^\eps_{k_M}$ is not reached by optimal
trajectories in the interval $]0,\pi]$. At the point
$C^\eps_{k_M}(0)$, an overlap
curve starts and passes through the origin (Fig.~\ref{f-TUTTALASINTESI}, Case B) .

\i{}\underline{CASE C:  $\ra\in(0,\frac{\al_1}{2\al})$}. The situation is more
complicated and it is depicted in the bottom of
Fig.~\ref{f-TUTTALASINTESI}, Case C.

\ed \ed {For $r=0$, the situation is the same as in CASE A, but for
the switching curve starting at $C^\eps_{k_M-1}(0)$.}

As explained in the introduction, the presence of
several cyclically alternating patterns of optimal synthesis, each of
them depending on an arithmetic property of $\al$, was already confirmed in \cite{q5},
by counting the number
of optimal trajectories reaching the south pole.

\brem
The first conjecture is implicitly disproved by the results of this paper.
More precisely an immediate consequence of our results is that the switching curve 
$C^\eps_{k_M-2}$ is always locally optimal, while $C^\eps_{k_M-1}$ is not, in general.
However, for every fixed $\bar r<\frac{1}{2}$ there exists $\al$ small enough with 
$\bar r\leq r(\al)<\frac12$ such that $C^\eps_{k_M-1}$ is locally optimal too, 
which contradicts the conjecture.
On the other hand Conjecture {\bf C2} is correct and, at the light of our main results, 
is completely proved and clarified.
\erem

\subsection{Notations}\label{not0}
All along the paper we use the notation $\eps=\pm1$. The set $so(3)$ of $3\times 3$
skew-symmetric matrices is a
three-dimensional vector space on which the following bilinear map
$$
<A,B>=-Tr(AB),\quad A,B\in so(3),
$$
is an inner product. For $A\in so(3)$, $\|A\|:=\sqrt{<A,A>}$ is the
norm (or length) of $A$. With the above notations, $F$ and $G$ are perpendicular
and normalized so that $\|F\|=\cos(\al)$ and $\|G\|=\sin(\al)$.

Let $Id$ be the $3\times 3$ identity matrix. We recall that $N=(0,0,1)^T$
and denote the south pole as $S=(0,0,-1)^T$. Set $c_t:=\cos(t)$
and $s_t:=\sin(t)$ for $t\in [0,2\pi)$. Recall that $X_+:=F+G$ and
$X_-:=F-G$ and we have
$$
X_+=\begin{pmatrix}
0 & -\ca & 0\\
\ca & 0 & -\sa\\
0 & \sa & 0
\end{pmatrix},\quad
X_-=\begin{pmatrix}
0 & -\ca & 0\\
\ca & 0 & \sa\\
0 & -\sa & 0
\end{pmatrix}.
$$
Let $\Pi_{x_3}$ be the orthogonal symmetry with respect to the $x_3$-axis,
i.e. $\Pi_{x_3}$ is represented in the canonical basis by
$Diag(-1,-1,1)$. Then, we have the following trivial but useful property.

\beq\label{sym0} \Pi_{x_3} X_\eps=X_{-\eps}\Pi_{x_3}. \eeq

We next recall standard formulas for a rotation $e^{tY}$ of $SO(3)$
in terms of its axis $Y$ (whose length is equal to one) and its
angle $t$. We have \beq\label{rot1} e^{tY}=
Id+s_tY+(1-c_t)Y^2. \eeq Moreover, for $t\in [0,2\pi)$, we
have
$$
e^{\Theta(t) Z_-(t)}:=e^{tX_+}e^{tX_-}\quad  e^{\Theta(t)
  Z_+(t)}:=e^{tX_-}e^{tX_+},
$$
where the unit vectors $Z_+(t),Z_-(t)$ are defined by
\beq\label{zpm}
Z_+(t)=\begin{pmatrix}
0 & -C(t) & -B(t)\\
C(t) & 0 & 0\\
B(t) & 0 & 0
\end{pmatrix},\quad
Z_-(t)=\begin{pmatrix}
0 & -C(t) & B(t)\\
C(t) & 0 & 0\\
-B(t) & 0 & 0
\end{pmatrix},
\eeq
with $B(t):=\frac{\sa s_{t/2}}{\sqrt{s^2_{t/2}\sa^2+c^2_{t/2}}}$,
$C(t):=-\frac{c_{t/2}}{\sqrt{s^2_{t/2}\sa^2+c^2_{t/2}}}$
and the angle
$\Theta(t)$ by
\beq\label{th0}
\Theta(t)=2 \arccos (s^2_{t/2}c_{2\al}-c^2_{t/2}).
\eeq

\section{The Extremal Front}
\subsection{Definition and description}\label{EXTRE}
As said in the introduction,
$\F(\al,T)$ the extremal front along $(\Sigma)_\al$
at time $T$ is the set of points reached at time $T$ by extremal
trajectories starting from $N$, i.e.
\bqn
\F(\al,T):=\{\bar x\in S^2:~~ \exists \mbox{ an extremal pair }
(x(.),\la(.))\mbox{ such that } x(0)=N,~~~x(T)=\bar x        \}.
\eqnl{eq-front}
Such extremals are parametrized by the length of the
first bang arc, the one of the last bang arc and the number of arcs:
\bqn
\Xi^+(s,t)&=&\overbrace{e^{X_\eps t}e^{X_{-\eps} v(s)}\cdots
e^{X_- v(s)}e^{X_+ s}}^{n ~~\mbox{terms}}N,\llabel{setEX-1}\\
\Xi^-(s,t)&=&\underbrace{e^{X_{\eps'} t}e^{X_{-\eps'} v(s)}\cdots
e^{X_+ v(s)}e^{X_- s}}_{n'~~\mbox{terms}}N,\llabel{setEX-2}
\eqn
where $s\in (0,\pi]$, $t\in (0,v(s)]$, the number of bang arcs  ($n$ and
$n'$ respectively) is an integer and
\bd
\i[(-)] $\eps=+1$ (resp.   $\eps=-1$), if $n$ is odd
(resp. even),
\i[(-)] $\eps'=+1$ (resp.   $\eps'=-1$), if  $n'$ is even
(resp. odd).
\ed
Roughly speaking, we would like to compute the limit, as $\al\to0$, of
$\F(\al,T)$, when $T$ is
such that the extremal front reaches
a \neigh of the south pole.

The idea is that, once one knows the extremal front $\F(\al,T)$ and if it is
optimal, then one can continue to build the synthesis for
times bigger than $T$ using  $\F(\al,T)$  as a \und{source} for the
minimization problem.

The identification of the front  $\F(\al,T)$ is not easy since it requires
the computation of the product of several exponentials of matrices. Moreover, if $\F(\al,T)$
crosses some switching curve, then the number of exponentials in
general depend on the point.

This problem is overcome by considering $\F(\al,T)$ only at times equal to
multiples of $\pi$. Indeed, first notice that, for
$T=\pi\left[\frac\pi{2\al}\right]$, the extremal front
reaches the points $C^\pm_{k_M}(0)$, i.e. the points where the last
switching curves $C^\pm_{k_M}$ start. Thanks to Proposition~\ref{ppp} below, at these
times, every extremal trajectory has the  same number of
switchings. The extremal front at times that are not multiple of $\pi$ can be
obtained {\sl a posteriori}, continuing the extremal front, as explained above.

From the structure of the extremal trajectories it follows that 
the time at which the point  $C^\pm_k(s)$ is reached is $T_k(s)=s+k v(s)$.
\bl
Let $k$ be an integer satisfying  $\displaystyle 1\leq k\leq\NMON:=\left[
\frac{({\cot (\alpha )}^2-1)^2}{2{\cot (\alpha )}^2-1}
\right]$, then $T_k(s)$ is a strictly increasing function of
$s$.
\el
\proof
It holds
\bqn
\dert{s}T_k(s)=\frac{1 + 2\,c_s\,{\cot (\alpha )}^2 + {\cot
(\alpha )}^4 +
  k\,\left( 2 + 2\,c_s\,{\cot (\alpha )}^2 \right)}
  {1 + 2\,c_s\,{\cot (\alpha )}^2 + {\cot (\alpha )}^4}.
\eqnl{eq-riccardomagro}
It is clear that the denominator of the above fraction is never vanishing on $[0,\pi]$ 
if $\al<\pi/4$.
On the other hand the numerator, as a function of $s$, reaches its minimum at $s=\pi$,
where it is equal to $({\cot (\alpha )}^2-1)^2-k(2{\cot (\alpha )}^2-1)$, and then
the conclusion follows easily.
\eproof

As a consequence, we obtain the following important corollary.
\bc
Let $k$ be an integer satisfying  $\displaystyle 1\leq k\leq\NMON$. If an
extremal trajectory is switching at time $T=k \pi$, then the length $s$ of
the first bang arc satisfies $s=\pi$.
\ec
Since for $\al$ small $k_M\leq\NMON$, then, for
$T=k\pi$ where $k$ is a positive integer such that $k\leq
[\pi/(2\al)]$, we have that all the extremal trajectories switch exactly
$k$ times (except the trajectories with length of the
first switching equal to $\pi$ that switch $k-1$ times). Therefore,
the extremal front $\F(\al,k\pi)$ is described by the next proposition.

\bp\label{ppp} Let $k$ be a positive integer such that $\displaystyle 1\leq
k\leq [\pi/(2\al)]$. Then, if $\al$ is small enough, we have

\bqn\label{uni} \F(\al,k\pi)=
\left\{\mathcal{E}^+(\al,k,s),~~~s\in(0,\pi]\right\}\bigcup
\left\{\mathcal{E}^-(\al,k,s),~~~s\in(0,\pi]\right\},\mbox{~~~where}: \eqn
\bqn
\mathcal{E}^+(\al,k,s):=\left\{\ba{l} e^{(k\pi-(k-1) v(s)-s) X_-}
\,e^{\frac{k-1}2\Theta(v(s))Z_-(v(s))}
\,e^{sX_+}N\mbox{~~~~~for $k$ odd,}\\
e^{(k(\pi-v(s))-s ) X_+}e^{\frac{k}2\Theta(v(s))Z_-(v(s))}e^{sX_+}N\mbox{~~~~~~~~~~~~for $k$
even.} \ea \right.
\eqnl{unidue}
The expression for $\mathcal{E}^-$ is the same as the
expression for $\mathcal{E}^+$ after exchanging the subscripts
$+$ and $-$. As a consequence, $\mathcal{E}^{-\eps}=\Pi_{x_3}\mathcal{E}^\eps$,
where $\Pi_{x_3}$ is the orthogonal symmetry with respect to the
$x_3$-axis. \ep

\brem\label{r11} 
Notice that  $\mathcal{E}^\eps(\al,k,0)=\mathcal{E}^{-\eps}(\al,k,\pi), \quad \eps=\pm\,,$ so that
$\F$ is described by a continuous closed curve.
\erem

\subsection{Description of the extremal front $\F(\al,k_M\pi)$ and consequences}
\llabel{des100}
As sketched in the introduction, we must describe the optimal
synthesis on $S^2$ deprived of a neighborhood of the south pole. For
that purpose, we will provide the precise asymptotics of
$\F(\al,k_M\pi)$, as $\al$ tends to zero, and derive, from its
topological nature, the minimum time front at time $k_M\pi$. 

From now on, for simplicity, we drop the dependence of $\mathcal{E}^\eps$ on $k_M$, 
i.e. we set $\mathcal{E}^\eps(\al,s):=
\mathcal{E}^\eps(\al,k_M,s)$, and we assume that $k_M$ is odd.\\

In the following, it will be useful to think of $\al$ and $r$ as two independent 
variables. For this purpose, define
$$
\ba{l}
\vspace{0.1cm}
\psi(\al,r,s):=\left(\frac{\pi}{2\al}-r\right)(\pi-v(s))+v(s)-s\,, \\
\vspace{0.1cm}
\theta(\al,r,s):=\left(\frac{\pi}{4\al}-\frac{1+r}{2}\right)\Theta(v(s))\,, \\
\chi^\eps(\al,r,s):=e^{\psi(\al,r,s)X_{-\eps}}e^{\theta(\al,r,s)Z_{-\eps}(v(s))}e^{sX_\eps}N\,.
\ea
$$
It is clear from \r{unidue} that 
$$\mathcal{E}^\eps(\al,s)=\chi^\eps(\al,r(\al),s)\,.$$
The following result is the key point in order to describe the extremal front at time $k_M\pi$.
\bl\label{fun}
There exists $\al_0>0$ such that the function $\chi^\eps\,,\ \eps=\pm$, 
defined above, is real-analytic for $(r,s,\al)\in \R^2\times I$, where $I=(-\al_0,\al_0)$.
Moreover, it admits a convergent power series
\beq\label{psw0}
\chi^\eps(\al,r,s)=\sum_{l\geq 0}f_l^\eps(s,r)\al^l,
\eeq
where the
$f_l^\eps(s,r)$ are real-analytic functions of $(s,r)\in \R^2$, 
$2\pi$-periodic in $s$ (therefore they are bounded over $\R\times [0,1]$).

As a consequence, the extremal front $\F(\al,k_M\pi)$, which is a continuous
closed curve, is piecewise analytic with discontinuities at
$s=0,\pi$ for derivatives of order greater than or equal to one.
\el


\noindent {\bf Proof of Lemma~\ref{fun}.}

We will prove the proposition only for $\chi^+$. 
Since $\chi^+$ is $2\pi$-periodic in $s$ and $r$ enters in an affine way in $\psi$ and $\theta$,
the real issue of analyticity revolves around the variable $\al$. First of all, it is clear that
$v(s)$ is actually a real-analytic function for $(s,\al)\in \R\times I$, where $I=(-\al_0,\al_0)$ 
with $\al_0>0$ small enough. Therefore, one has only to prove
the real-analyticity of $\tilde\psi(\al,s):=\frac{v(s)-\pi}{\al}$ and 
$\frac{\beta(s,\al)}{\al}$, where $\beta(s,\al):=\Theta(v(s))$, for $(s,\al)\in \R\times I$, 
where $I=(-\al_0,\al_0)$, for some $\al_0>0$. 

Note that 
$$
\tilde\psi(\al,s)=\frac2{\al}\arctan(\sac\mu(s)),\ \ \
\mu(s):= \frac{\sis}{\cac+\sac\cs}.
$$
The function $\mu$ is real-analytic for $(s,\al)\in \R\times I$, with $I$ open neighborhood of zero,
and thus uniformly bounded over $\R\times I$. In addition, $\arctan(\cdot)$ is real analytic in a 
neighborhood of zero. Hence the conclusion for $\tilde\psi(\al,s)$.

As for $\frac{\beta(s,\al)}{\al}$, rewrite first Eq.~\r{th0} as
$$
\cos(\beta(s))=1-G(s,\al),
$$
with
\beq\label{koko}
G(s,\al):=2\sac\left[1+\frac{\cac\sac\mu^2(s)}{1+\sa^4\mu^2(s)}+2\sac
(1+\frac{\cac\sac\mu^2(s)}{1+\sa^4\mu^2(s)})^2\right].
\eeq
We first need to determine a convergent power series for $\beta$ from
the expression 
\beq\label{eo0-2}
\beta=\arccos(1-G).
\eeq
Note that $|G(s,\al)|\leq
5\al^2$ for $\al$ small enough. We first expand $\arccos(1-G)$ in a power 
series in $G$. Starting from the power series
$$
(1-t)^{-1/2}=1+\sum_{m\geq 1}s_mt^m,
$$
with radius of convergence equal to $1$ we get 
$$
\frac{d}{dG}\big(\arccos(1-G)\big)=-\frac1{\sqrt{2G}}\frac1{\sqrt{1-G/2}},
$$
and, after simple integration, 
\beq\label{eo0-1}
\arccos(1-G)=-\sqrt{2G}\big(1+\sum_{m\geq
1}\frac{s_m}{2^{m+1}(m+1/2)}G^m\big),
\eeq
Finally, from Eq.~\r{koko}, $G$ can be written as $2\sac\big(1+\sac H(s,\al)\big)$
with $H(s,\al)$ uniformly bounded by $3$. Then,
\beq\label{eo0}
\sqrt{2 G(s,\al)}= 2\sa\big(1+\sac H(s,\al)\big)^{1/2}.
\eeq
Gathering Eqs.~\r{eo0-2},\r{eo0-1},\r{eo0}, we get the real-analyticity of
$\frac{\beta(s,\al)}{\al}$ for $(s,\al)\in \R\times I$, 
where $I=(-\al_0,\al_0)$, for some $\al_0>0$ small enough. 

\eproof

We next compute $f_0,f_1,f_2$ and obtain the following proposition.
\bp\label{pop0} 
For $\al$ small enough, the function $\chi^\eps\,,\ \eps=\pm$ defined above
and its derivative with respect to $s$ have the following expansion
\bqn
\chi^\eps(\al,r,s)&=&f^\eps_0(s,r)+f^\eps_1(s,r)\al+f^\eps_2(s,r)\al^2+\O(\al^3)\,,\llabel{ser}\\
\derp{s}\chi^\eps(\al,r,s)&=&\derp{s}f^\eps_0(s,r)+\derp{s}f^\eps_1(s,r)\al+
\derp{s}f^\eps_2(s,r)\al^2+\O(\al^3)
\eqnl{derser}
where $f^\eps_l\,,~l=0,1,2$ are defined as in \r{psw2} and $|\O(\al^3)|\leq C|\al^3|$,
with the constant $C$ independent of $s\in\R$ and $r\in[0,1)$.
\ep

\noindent {\bf Proof of Proposition~\ref{pop0}.}
We will prove the proposition only for $\chi^+$. 
To proceed, we list the expansions of the form \r{ser} for several quantities, 
obtained after elementary computations.
\bqn
\psi(\al,r,s)&=&\pi-s-\pi\sis \al+2(1-r)\sis\al^2+\O(\al^3)\llabel{psi1}\\
\theta(\al,r,s)&=&\pi-2\al(1+r)+\frac{\pi\sisc}2\al^2+\O(\al^3)\llabel{kthe1}\\
Z_-(v(s))&=&\bpm 0&-\al\sis&1-\frac{\sisc}2\al^2\\\al\sis&0&0\\
-1+\frac{\sisc}2\al^2&0&0\epm+\O(\al^3)\llabel{Z-1}.
\eqn
Using Eqs.~\r{psi1},\r{kthe1}, we get that
\bqn
\sin\big(\psi(\al,r,s)\big)&=&\sis+\pi\sis\cs\al-\big(2(1-r)\sis\cs+
\frac{\pi}2 s_s^3\big)\al^2+\O(\al^3)\llabel{spsi}\\
\cos\big(\psi(\al,r,s)\big)&=&\cs-\pi\sisc\al+\big(2(1-r)\sisc-
\frac{\pi^2}{2}\sisc\cs\big)\al^2+\O(\al^3)\llabel{cpsi}\\
\sin\big(\theta(\al,r,s)\big)&=&2\al(1+r)-\frac{\pi\sisc}2\al^2+\O(\al^3)\llabel{skthe1}\\
\cos\big(\theta(\al,r,s)\big)&=&-1+2\al^2(1+r)^2+\O(\al^3)\llabel{ckthe1}. \eqn
Using Eqs.~\r{rot1},~\r{spsi} and \r{cpsi} we obtain \beq\llabel{int1}
e^{\psi(\al,r,s)X_-}=\bpm
-\cs+\pi\sisc\al&-\sis-\pi\cs\sis\al&-(1+\cs)\al\\
\sis+\pi\sis\cs\al&\cs+\pi\sisc\al&\sis\al\\
-(1+\cs)\al&-\sis\al&1 \epm+\mathcal{R}(s)\al^2+\O(\al^3)\,,
\eeq
where
$$
\mathcal{R}(s)=\bpm
\cs+\frac{\pi^2}2\cs\sisc+1-2(1+r)\sisc & \frac\sis 2+2\cs\sis (1+r)+\frac{\pi^2}2 s^3_s & \pi \sisc \\
-\frac\sis 2-2\cs\sis (1+r)-\frac{\pi^2}2 s^3_s & -2 (1+r)\sisc+\frac{\pi^2}2\cs\sisc &  \pi\sis\cs\\
\pi\sisc & -\pi\sis\cs & -1-\cs
\epm\,,
$$
and using Eq.~\r{rot1}, we have \beq\llabel{I1} e^{sX_+}N=\bpm
\sa\ca(1-\cs)\\-\sa\sis\\1-\sa^2(1-\cs)\epm= \bpm
\al(1-\cs)\\-\al\sis\\1-\al^2(1-\cs)\epm+\O(\al^3). \eeq
An easy computation yields \beq\llabel{Z2-1} Z_-^2(v(s))=\bpm
-1&0&0\\0&-c^2(s)&b(s)c(s)\\0&b(s)c(s)&-b^2(s)\epm= \bpm
-1&0&0\\0&-\al^2\sis^2&\al\sis\\0&\al\sis&-1+\al^2\sis^2\epm+\O(\al^3).
\eeq
Using Eqs.~\r{rot1}, \r{skthe1} \r{ckthe1} and \r{I1} and the previous equation, we get
\bqn
\label{I2}
e^{\theta(\al,r,s)Z_-(v(s))}e^{sX_+}N= \left(\ba{c}
\alpha(1  + 2\,r   + c_s) - \al^2 \frac{\pi}{2}\,s_s^2\\
\alpha \,s_s\\
-1 + {\alpha }^2\,\left( 1 + 2\,r + 2\,r^2 + c_s + 2\,r\,c_s
\right)
\ea\right)+\O(\al^3).
\eqn
Applying $e^{\psi(s)X_-}$ to
the previous equation and using Eq.~\r{int1}, we
finally get Eq.~\r{ser} for $\eps=+$. The expression of the derivative~\r{derser} is then an immediate 
consequence of the analiticity of the function $\chi^+$.

The results for $\chi^-$ and $\derp{s}{\chi^-}$
are obtained similarly together with Eq.~\r{sym0}.



\eproof
 
Since the quantities of the form $\O(\al^3)$ in Proposition~\ref{pop0} satisfy
$|\O(\al^3)|\leq C|\al^3|$ for some $C$ independent of $r$, the expressions 
\r{yacineseiunbalordo}--\r{derfronte} are straightforward consequences. 
Hence the shape of the extremal front at time $T=k_M\pi$ is known for $\al$ 
small. In particular its image with respect to the map $M_{\al}$ defined in 
Section~\ref{intro} is approximated, in the $\con^1$ sense, by a circle of 
radius $2\,r(\al)$ centered at the origin.

We finally note that, for $k_M$ even, Lemma~\ref{fun} is still valid, while, 
with computations similar
to those made in the proof of Proposition~\ref{pop0}, it is easy to see that the 
formulas for $f_k^\eps,\ k=0,1,2$ simply differ, with respect to \r{psw2}, for the 
sign of the first two components.

%
%
%

\brem\label{lem-est} Repeating the previous computations, we also obtain series
expansions for $\mathcal{E}^\eps(s,k_M-1,\al)$ and $\derp{s}\mathcal{E}^\eps(s,k_M-1,\al)$. 
Indeed, we just have to replace $r$ by $1+r$.
In that case the shape of the extremal front $\F(\al, (k_M-1)\pi)$, after applying 
the map $M_{\al}$, is approximated, in the $\con^1$ sense, by a circle of radius 
$2(1+r(\al))$ centered at the origin.
\erem

\section{Case $r(\al)=\bar{r}\in (0,1)$}\label{barr}
In this section, we study the case in which $\al$ tends to zero with $r(\al)=\bar{r}$,
for a constant $\bar{r}\in (0,1)$. More precisely we consider the decreasing sequence 
$\al_k=\frac{\pi}{2(k+\bar{r})}$, for $k\geq 1$. We first
describe the minimum time front at $T=k_M\pi$, then we identify and study
the candidate for the limit synthesis and finally we prove Theorem~\ref{the1}.

\subsection{Description of the minimum time front at $T=k_M\pi$}\llabel{ss-des}
The purpose of the paragraph is to prove the following proposition.
\bp\label{pop2}
Fix $\de>0$ small. For $\al$ small enough with $r(\al)>\de$, the 
extremal front $\F(\al,k_M\pi)$ is homeomorphic to a circle. 
As a consequence, the switching curves defined inductively in 
Eq.~\r{swi0} are optimal up to $k=k_M$ and $OF(\al,k_M\pi)$, the 
minimum time front at time $k_M\pi$ coincides with $\F(\al,k_M\pi)$.
\ep

\noindent {\bf Proof of Proposition~\ref{pop2}.}
From Proposition~\ref{pop0} 
we get that
the extremal front $\F(\al,k_M\pi)$
is the union of two arcs, $\mathcal{E}^+(\al,s)$, $s\in[0,\pi]$ and
$\mathcal{E}^-(\al,s)$, $s\in[0,\pi]$ so that, for $\eps=\pm$ and $s\in [0,\pi]$,
\beq\label{f+0}
\mathcal{E}^\eps(\al,s)=\bpm -2r(\al)\eps \al\cs\\ 2r(\al)\eps \al\sis\\-1\epm
+\O(\al^2),
\eeq
and
\beq\label{df+0} \derp{s}\mathcal{E}^\eps(\al,s)=2r(\al)\eps\al\bpm \sis\\
\cs\\0\epm +\O(\al^2).
\eeq
Moreover, at $s=0$ and $s=\pi$, the derivatives of $\mathcal{E}^\eps(\al,s)$
are only one-sided, i.e. as $s>0$ tends to zero and $s<\pi$ tends to
$\pi$. By a trivial continuity argument,
one can parameterize $\F(\al,k_M\pi)$ as a closed continuous curve
$\gamma$ defined on $[0,2\pi]$ so that $\gamma(s)=\mathcal{E}^+(\al,s)$
for $s\in(0,\pi]$ and $\gamma(s)=\mathcal{E}^-(\al,s-\pi)$ for
$s\in(\pi,2\pi]$. Moreover, with the previous computations, it is
immediate that $\gamma$ is in fact piecewise $C^1$ with possible
discontinuity jumps for $\dert{s}\gamma$ at $s=0$ and $s=\pi$.

Since the curve $\gamma$ is in a neighborhood of the south pole of
size proportional to $\al$ (thanks to Eq.~\r{f+0}),
it is enough to prove that the orthogonal projection $\gamma_1$ of
$\gamma$ on the $(x_1,x_2)$-plane is homeomorphic to the circle
$e^{is}$, $s\in [0,2\pi]$. Using Eq.~\r{f+0}, we see that
$\|\gamma_1(s)\|=2r(\al)\al+\O(\al^2)$ on $[0,2\pi]$, which implies that
the continuous function $\|\gamma_1(s)\|$ is always strictly
positive for $\al$ small enough. We can therefore parameterize
$\gamma_1$ using polar coordinates $(\rho,\beta)$, i.e., for $s\in
[0,2\pi]$,
$$
\gamma_1(s)=\rho(s)e^{i\beta(s)},
$$
where $\rho(\cdot):=\|\gamma_1(\cdot)\|$ and the function $\beta(\cdot)$ are
defined on $[0,2\pi]$, continuous and piecewise $C^1$, with possible
jumps of discontinuity for their derivatives at $s=0$ and $s=\pi$.\\
In addition $\rho(0)=\rho(2\pi)$, $\beta(0)\equiv\beta(2\pi)\equiv\pi$ (mod $2\pi$) and,
from Eq. \r{f+0}, $\beta(s)=\pi-s+\O(\al)$.
To prove Proposition~\ref{pop2}, it suffices now to
prove that $\beta$ is a monotone bijection from $[0,2\pi]$ to $[-\pi,\pi]$. 
The latter simply results from
Eq.~\r{df+0}. Indeed, from that equation, we get that
$\dert{s}\beta(s)=-1+\O(\al)$ where $\beta$ is differentiable and
the one-sided derivatives at $s=0$ and $s=\pi$ verify the same
equation. We deduce that $\beta$ is strictly decreasing for $\al$ small enough.

We next show that $OF(\al,k_M\pi)$, the minimum time front at time
$k_M\pi$ coincides with $\F(\al,k_M\pi)$. By the results of
\cite{q5}, we first notice that any time minimal trajectory starting
at the north pole reaches the south pole in time $T>k_M\pi$.
Therefore $OF(\al,k_M\pi)$ is not empty and is included in
$\F(\al,k_M\pi)$ according to the PMP. According to Theorem $27$ of
\cite{libro}, $OF(\al,k_M\pi)$ is a one-dimensional piecewise $C^1$
compact embedded submanifold of $S^2$. By an easy topological
argument, we deduce from the above that
$OF(\al,k_M\pi)$ coincides with $\F(\al,k_M\pi)$. \eproof

\brem\label{opt0}
Thanks to Remark~\ref{lem-est}, and with arguments similar to those 
of the previous proof, one can prove that $\F(\al, (k_M-1)\pi)$ is optimal
for $\al$ small enough, with no assumptions on the remainder $r$.
\erem

\subsection{Optimal synthesis for the linear pendulum control problem}\llabel{s-rwlp}
Recall that $M_\al:\R^3\rightarrow \R^2$ is the composition of the projection
$(x_1,x_2,x_3)\mapsto (x_1,x_2)$
followed by the dilation by $1/\al$. With the results of the
previous subsection, it is clear that the original control problem on
$S^2$ can be reduced, near the south pole, to a planar
control problem on the neighborhood of the south pole delimited by
$\widetilde{OF}(\al,k_M\pi):=M_\al(OF(\al,k_M\pi))$ along $(\widetilde{\Sigma})_{\al}$, 
the control system obtained as the image
of $(\Sigma)_\al$ by $M_\al$, i.e.
\beq\label{f-al} (\widetilde{\Sigma})_{\al}:\quad\left\{ \ba{l}
  \dot{z_1}=-\cos(\al)  z_2, \\
  \dot{z_2}=\cos(\al) z_1 +u\frac{\sin(\al)}\al \sqrt{1-(\al z_1)^2-(\al z_2)^2},
\ \ (z_1,z_2)\in\R^2,\ \ |u|\leq 1.
\ea \right.
\eeq
It is therefore natural to conjecture (simply set $\al=0$ in $\widetilde{OF}(\al,k_M\pi)$ and
$(\widetilde{\Sigma})_\al$) that the limit synthesis should be that of connecting 
the circle of radius $2r(\al)$, $C(0,2\ra)$,
to every point of the disk $B(0,2\ra)$ along the control system $(Pen)$ given by
Eq.~\r{lineare0}, which we rewrite as
\bqn
(Pen)\quad \dot{z}=A_0z+ub_0,\quad\mbox{ with }\quad
A_0=\bpm 0&-1\\1&0\epm,\quad b_0=\bpm 0\\1\epm,
\eqnl{lineare}
where $z\in\R^2$ and $u\in [-1,1]$. The control system $(Pen)$ corresponds to
a linear pendulum with a forcing term.

Theorem~\ref{the1} simply states that the conjecture is correct and, as a first step
for an argument, we describe, in more details in this subsection, the conjectured limit
synthesis. 
Hence we focus on the following problem.\\

{\bf (P)} {\sl Fixed $\rho\in]0,2]$, for any given $\bar{y}\in B(0,\rho)$
find a time optimal trajectory connecting the circle of radius $\rho$ centered at the
origin to $\bar{y}$ along the control system $(Pen)$.}

\brem 
The problem of computing the TOS for the linear pendulum, taking the origin as a source, is studied
in any textbook of optimal control. Here as source we take 
the circle of radius $\rho$ centered at the origin, which is a level set of the Hamiltonian 
$H=\frac12(z_1^2+z_2^2)$ associated to the uncontrolled system.
\erem 
It is easy to see that the solutions of problem {\bf (P)} must be bang-bang trajectories. 
Indeed since $(Pen)$ is a bidimensional linear control system it is well known that this 
property is guaranteed by the Kalman controllability condition $\det(b_0,A_0b_0)\neq 0$, which is satisfied
by $(Pen)$.
To determine the TOS,
we first look for the switching curves. We know that every extremal
trajectory for the problem {\bf (P)} must satisfy the transversality
condition of the PMP stated in Section~\ref{ss-basic}.
Here the source manifold is the circle $C(0,\rho)$, and the transversality 
condition essentially translates into the property that the vector 
$\lambda(0)$ (that, without loss of generality, we will assume unitary) 
is proportional to $z(0)\in C(0,\rho)$ (identifying the cotangent space 
with the plane $\R^2$). To determine completely $\lambda(0)$, it is enough to
observe that a necessary condition for $z(.)$ to be optimal is that
$\dot{z}(0)$ points inside the disk $B(0,\rho)$, i.e., if we denote by
$u_{opt}$ the corresponding control, then
$$<z(0),\dot{z}(0)>\leq 0\iff
<z(0),A_0z(0)+u_{opt}b_0>\leq 0\iff
<z(0),u_{opt}b_0>\leq 0\,.
$$
Therefore,
$u_{opt}=-$sgn$<z(0),b_0>$. On the other hand, from the
maximality condition of the PMP, we have
$u_{opt}=$sgn$<\lambda(0),b_0>$ and, therefore, one can define
$\lambda(0):=-z(0)/\rho$.  Finally $u_{opt}=-$sgn($z_2(0)$) (except at
the points $\pm(\rho,0)$ ), while the switching time $t_{sw}$
must satisfy the condition  $<\lambda(t_{sw}),b_0>=\lambda_2(t_{sw})=0$.\\
Consider now the adjoint system \bqn \left\{ \ba{l}
  \dot{\lambda}_1=-\lambda_2, \\
  \dot{\lambda}_2=\lambda_1. \ea \right.  \eqnl{adj} If we identify
$\R^2$ with the complex plane, so that $z=z_1+iz_2$ and
$\lambda=\lambda_1+i\lambda_2$, then the equations \r{lineare},
\r{adj} become $$\dot{z}=i(z+u)\qquad\mbox{ and
}\qquad\dot{\lambda}=i\lambda.$$
Moreover we can set $z(0)=-\rho
e^{-i\theta}$ and $\lambda(0)=e^{-i\theta}$ for some
$\theta\in[0,2\pi[$ and the corresponding solutions are:
$$
\left\{
\ba{l}
z(t)=(z(0)+u_{opt})e^{it}-u_{opt}=-\rho e^{i(t-\theta)}+u_{opt}(e^{it}-1),\nn\\
\lambda(t)=\lambda(0)e^{it}=e^{i(t-\theta)}.\nn
\ea
\right.
$$
The switching curves are determined by the relation
$t_{sw}\equiv\theta$ (mod $\pi$) and this allows to conclude that
the switching curves are the following two semicircles of radius 1:
\bqn \left\{ \ba{l}
  z(\theta)=1-\rho-e^{i\theta}\qquad \theta\in [0,\pi[,\\
  z(\theta)=\rho-1-e^{i\theta}\qquad \theta\in [\pi,2\pi[. \ea
\right.\nn \eqn
These switching curves cannot
be optimal for $\rho<2$ since they are not
locally optimal, as can be easily checked using the definition given in 
Section~\ref{des000}. 
We conclude that the optimal trajectories are bang arcs and
the corresponding control depends on the sign of the component
$z_2(0)$ of the starting point.

To conclude the description of the synthesis, it is enough to determine
the cut locus, i.e. the set of points that are reached by two or more
optimal trajectories at the same time.  Assume that $z\in\C$ belongs
to the cut locus. Then, there exist
$s\in [0,\pi),s'\in [\pi,2\pi)$ and $t$ such that 
\bqn \left\{ \ba{l}
  z=-\rho e^{i(t-s)}+1-e^{it},  \\
  z=-\rho e^{i(t-s')}-1+e^{it}.  \ea 
\right.  \eqnl{overlap0}
Therefore $|z-1+e^{it}|=|z+1-e^{it}|=\rho$. In particular, denoting by
$\bar{z}$ the complex conjugate to $z$, we have \bqn
(z-1+e^{it})(\bar{z}-1+e^{-it})-(z+1-e^{it})(\bar{z}+1-e^{-it})=-4z_1+4z_1\cos
t+4z_2\sin t=0,
\label{overlap1}\\
(z-1+e^{it})(\bar{z}-1+e^{-it})+(z+1-e^{it})(\bar{z}+1-e^{-it})=
2z_1^2+2z_2^2+4-4\cos t=2\rho^2.
\label{overlap2}
\eqn From \r{overlap1} we have that $\cos
t=\frac{z_1^2-z_2^2}{z_1^2+z_2^2}$, and, substituting in \r{overlap2},
we find that $z$ must satisfy the equation
\bqn
z_1^4+z_2^4+2z_1^2z_2^2-\rho^2z_1^2+(4-\rho^2)z_2^2=0.
\eqnl{over-locus}
The previous computation show that the cut locus is a subset of the set of
points belonging to the locus defined by \r{over-locus}.
Actually it is easy to see that this is the proper subset obtained 
with the additional condition $z_1z_2\geq 0$, that corresponds to $t\leq \pi$.   
The precise shape of
the optimal synthesis, which is now clear, is portrayed in Fig.~\ref{sintesilineare} for
a particular value of $\rho<2$.
Notice that, from the previous computations, we have $\rho
e^{is'}=\rho e^{is}+2-2e^{it}$ and, since $\rho
e^{is'}+\rho e^{is}=2\rho e^{is'}-2+2e^{it}$ and
$\rho e^{is'}-\rho e^{is}=2-2e^{it}$ are orthogonal in
the complex plane, we find easily the following equation:
$$(2-\rho\cos s')(\cos t-1)-\rho\sin s'\sin t=0.$$
Consequently, for $t\in[0,2\pi[$ and $s'\in[\pi,2\pi[$, one has,
along the overlap curve \bqn t=t(s')=-2\arctan\frac{\rho\sin
s'}{2-\rho\cos s'}. \eqnl{over-param} This expression will be useful
in the following. Also, notice that combining \r{overlap0} and
\r{over-param} one easily finds a parametrization of the
overlap curve in terms of $s'$ and that in an analogous way it
is possible to parameterize it by means of $s$.
From now on we will denote by $\gamma^o_{pen}(\cdot)$ the parameterization 
of the overlap curve with respect to the parameter $s$.
\brem
If $\rho=2$  the previous reasoning 
does not apply and indeed the synthesis is different.
In this case the overlap curve coincides with the switching
curves and with the trajectories reaching the origin corresponding to $u=\pm 1$. 
A simple way to prove this fact is to study the optimal synthesis starting from the origin
with vector fields with opposite signs, and observe that the extremal front at time $\pi$
is a circle of radius $2$.
\erem

\ppotR{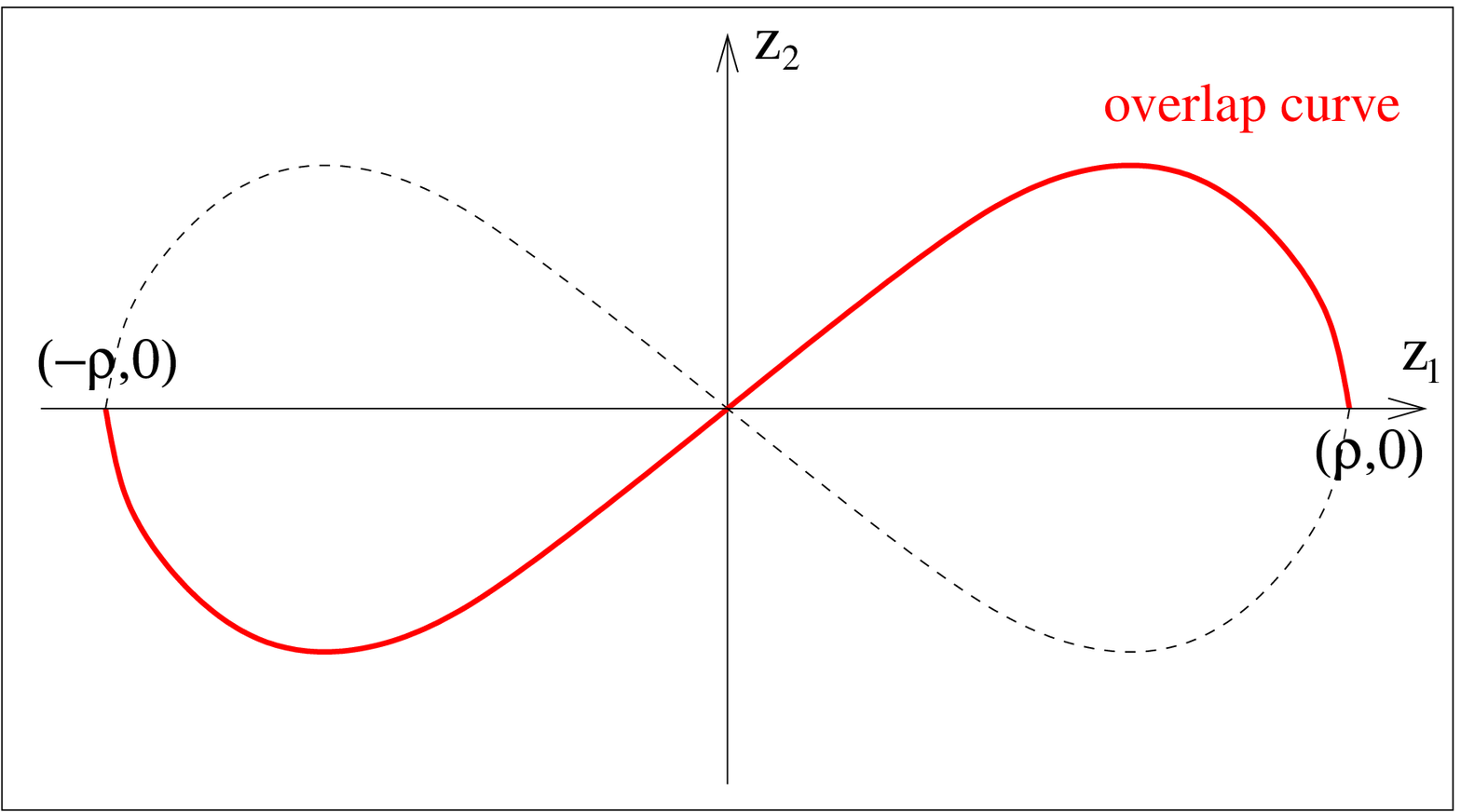}{The overlap curve for the pendulum problem}{8}

\subsection{Proof of Theorem~\ref{the1}}
The proof of Theorem~\ref{the1} is divided in two parts. Roughly speaking, 
defining $P_\eps^{\al}:=M_{\al}(C^\eps_{k_M}(0))$ for $\eps=\pm$, we will look separately at  
the shape of the synthesis far from $P_\eps^{\al}$ and inside neighborhoods of 
$P_\eps^{\al},\ \eps=\pm$. Let us call $V_\al$ the image with respect to $M_\al$ of the 
neighborhood of the south pole enclosed by $OF(\al,k_M\pi)$ and $B_\eps^{\al}(\xi)$, the ball of 
center $P_\eps^{\al}$ and radius $\xi$. 

Then, the previous cases correspond to the following two propositions whose meaning is clarified 
by Figure~\ref{f-proof}.

\ppotR{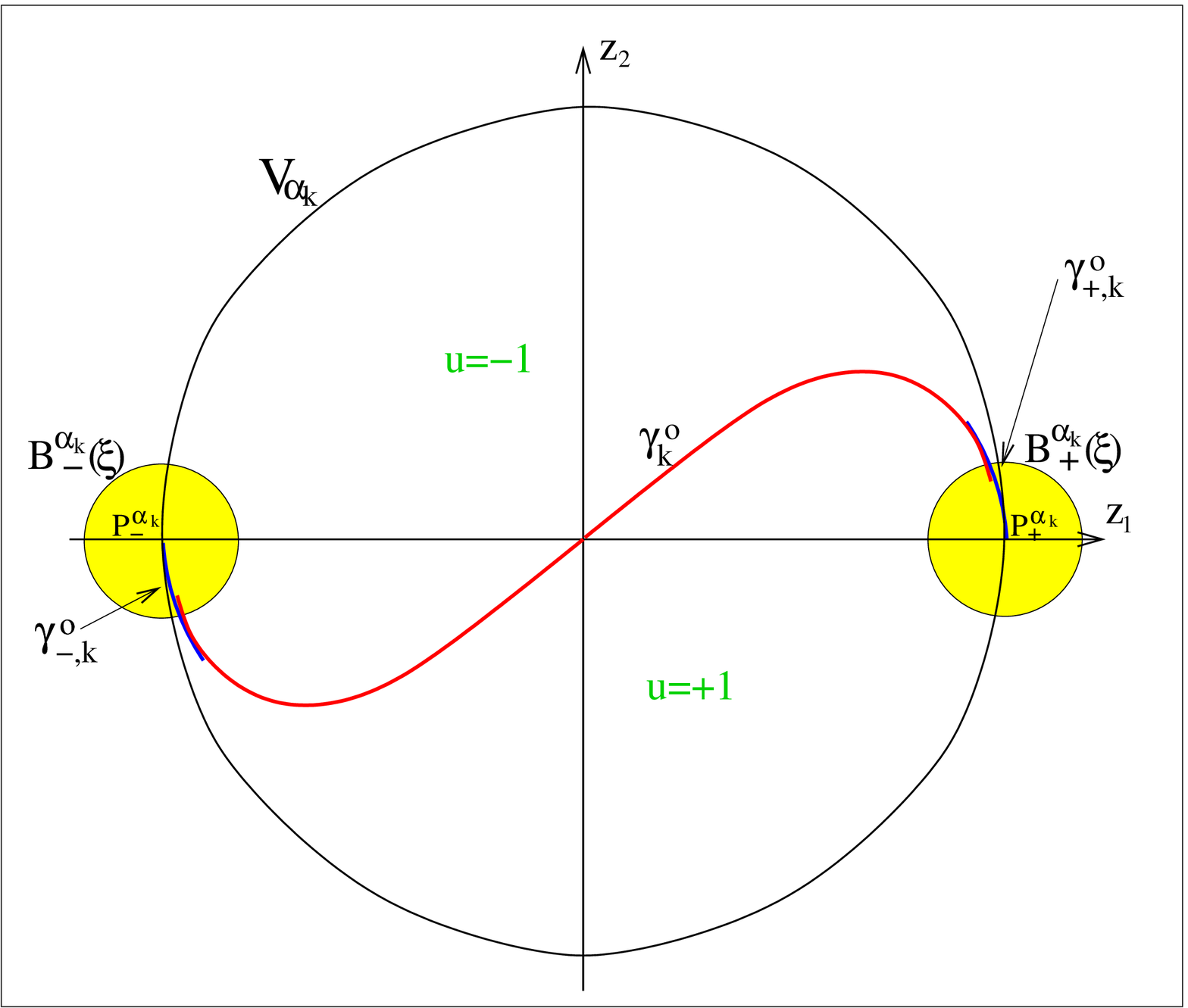}{Propositions \ref{p1-th} and \ref{p2-th}}{8}

\bp\label{p1-th}
Let $\bar{r}\in (0,1)$ and  $\al_k:=\frac{\pi}{2 (k+\bar r)}$. Then, for any $\xi>0$ there exist
a positive integer $\bar k$ and a compact interval $I\subset (0,\pi)$ such that it is possible 
to find a curve $\gamma^o_k$, defined on $I$ for $k\geq\bar k$, verifying the following: 
$\gamma^o_k$ divides 
$V_{\al_k}\setminus \big(B_+^{\al_k}(\xi)\cup B_-^{\al_k}(\xi)\big)$ in two connected 
components, such that ``above'' $\gamma^o_k$  
the optimal feedback associated to the synthesis for $\al=\al_k$ takes the value $-1$ and, below 
$\gamma^o_k$, it is equal to $1$, and in particular $\gamma^o_k$ is an overlap curve for $\al=\al_k$.
Moreover, $\gamma^o_k$ converges to $\gamma^o_{pen}$ in the
$C^0$ topology of $I$.
\ep

\bp\label{p2-th}
Consider the notations defined above. Then there exist $\xi>0$, $\tau_\eps,\ \eps=\pm$ with
$0<\tau_-<\tau_+<\pi$ and a positive integer $\bar k$ such that, for every $k\geq \bar k$, it is possible to find two curves 
$\gamma^o_{-,k}$ and $\gamma^o_{+,k}$, defined respectively on $[0,\tau_-]$ and
$[\tau_+,\pi]$, verifying the following:
$\gamma^o_{\eps,k}$ divides $V_{\al_k}\cap B_\eps^{\al_k}(\xi)$ in two connected components, such that ``above'' $\gamma^o_{\eps,k}$  
the optimal feedback associated to the synthesis for $\al=\al_k$ takes the value $-1$ and, below 
$\gamma^o_{\eps,k}$, it is equal to $1$, and in particular the $\gamma^o_{\eps,k}$ are overlap curves for $\al=\al_k$. Moreover, 
$\gamma^o_{-,k}$ and $\gamma^o_{+,k}$ converge to $\gamma^o_{pen}$ in the
$C^0$ topology respectively of $[0,\tau_-]$ and $[\tau_+,\pi]$.
\ep

The choice of studying the synthesis separately in neighborhoods of $P^\al_\eps$ and far
from $P^\al_\eps$ is justified by the fact that the proofs of the previous propositions rely on 
different implicit function arguments.

It is clear that, combining Proposition~\ref{p1-th}, for an appropriate choice of $\xi$,
with Proposition~\ref{p2-th}, one almost completes the proof of Theorem~\ref{the1}. 
We will not prove explicitly that the convergence of $\gamma^o_k$ to $\gamma^o_{pen}$
with respect to the parameter $\bar r$ is uniform in any closed interval of $(0,1)$. As explained in
Remark~\ref{extend}, this can be done with the same methods used in the proofs of
Propositions \ref{p1-th},\ref{p2-th}, but with much more computational efforts.

We will therefore only provide the complete proofs of the propositions.
For this purpose, the first step consists in checking whether
the switching curves $C^\eps_{k_M},\ \eps=\pm 1$ are optimal or not.
In that regard and similarly to the case of the linear pendulum, we have the
following result:
\bl
Let $\bar r\in ]0,1[$. Then, if $\al$ is small enough and $r(\al)=\bar r$,
the switching curve $C^\eps_{k_M}$ is nowhere locally optimal i.e. all the extremal trajectories 
switching on $C^\eps_{k_M}$ lose optimality before reaching it.
\llabel{lemma-nonlocop}
\el

\noindent {\bf Proof of Lemma~\ref{lemma-nonlocop}.}
For simplicity we define $S(s):=C^+_{k_M}(s)$ and we assume $k_M$ odd.
As in the proof of Lemma~\ref{lem-est}, we get the following asymptotic
expansions, after applying the map $M_\al$:
\bqn
S(s)&=&\left(\ba{c} 2r-1+c_s\\ s_s \ea\right)+\O(\al),\qquad\quad
S(0)=\left(\ba{c} 2r+\O(\al)\\ 0\ea\right),\\
S'(s)&=&\left(\ba{c} -s_s\\ c_s \ea\right)+\O(\al),\qquad \qquad\qquad
S'(0)=\left(\ba{c} 0\\ 1+\O(\al)\ea\right),\\
S''(s)&=&\left(\ba{c} -c_s\\ -s_s \ea\right)+\O(\al).
\eqn
Integrating the above equation, we have
\bqn
S'(s)&=&S'(0)+\int_0^s S''(\tau)d\tau=\left(\ba{c} -s_s+\O(s\al)\\ c_s+
\O(\al) \ea\right),\llabel{improve1}\\
S(s)&=&S(0)+\int_0^s S'(\tau)d\tau=\left(\ba{c} 2r-1+c_s+\O(\al)\\ s_s+
\O(s\al), \ea\right)\llabel{improve2}
\eqn
and therefore
$$\frac{1}{c_\al}X_\pm(S(s))=\left(\ba{c} -S_2(s)\\ S_1(s)\pm \frac{\tan\al}
{\al} \sqrt{1-\al^2S_1(s)^2-\al^2S_2(s)^2}\ea\right)=\left(\ba{c}
-s_s+ \O(s\al)\\ 2r-1+c_s+\O(\al)\pm(1+\O(\al^2))\ea\right).
$$ Here $S_i$, $i=1,2,$ denotes the $i$-th component of $S$. Dividing the
above equation by $1+\O(\al)$, we can assume that the first
component is identically equal to $-s_s$. The same can be done with
the expression \r{improve1}, so that it is possible to compare the
three vectors obtained in this way simply by looking at the second
components, which are equal respectively to $2r-1+c_s\pm
1+\O(\al)$ and $c_s+\O(\al)$. In particular, the fact that
$S(\cdot)$ is nowhere locally optimal if $\al$ is small enough
follows from the inequalities 
$2r-2+c_s+\O(\al)<c_s+\O(\al)
<2r+c_s+\O(\al).$

\eproof

A straightforward consequence of the previous result is the presence of a non
trivial cut locus in the neighborhood of the south pole enclosed by $F(\al,k_M\pi)$.
It remains to clearly define that cut-locus, which is the purpose of
Propositions~\ref{p1-th} and \ref{p2-th}. 

\subsubsection{Proof of Proposition~\ref{p1-th}}
As usual, we only provide an argument in the case $k_M$ odd and we fix
the remainder equal to $\bar r\in(0,1)$. 

Recall that, according to Section~\ref{ss-des}, $OF(\al,k_M\pi)$ is approximately (up to order $\al^2$) 
a circle of radius $2\bar r \al$. To describe the synthesis inside  the
neighborhood of the south pole enclosed by $OF(\al,k_M\pi)$, it is
more convenient to use the two dimensional control system
$(\widetilde{\Sigma})_{\al}$, which is rewritten as follows by using
Eq.~\r{f-al},
$$
\dot z=\ca A_0z+u\frac{\sa}{\al}\sqrt{1-\al^2\|z\|^2}b_0.
$$
We set $\tilde{X}^\al_\eps(z):=\ca A_0z+\eps
\frac{\sa}{\al}\sqrt{1-\al^2\|z\|^2}b_0$ and $\tilde{X}^{pen}_\eps(z):=A_0 z+\eps b_0$, for $\eps=\pm$,
and we define $\widetilde{OF}(\al,k_M\pi)$ as the image by $M_\al$ of $OF(\al,k_M\pi)$. Then we know that, up to order $\al$,   $\widetilde{OF}(\al,k_M\pi)$ is a circle of radius $2\bar r$.
In particular, as in the proof of Proposition~\ref{pop2}, one can construct a piecewise smooth
parameterization $\sigma_{\al}:[0,2\pi]\rightarrow
\Tilde{F}(\al,k_M\pi)$ so that $\sigma_{\al}(0)=P_-^{\al}$,
$\sigma_{\al}(\pi)=P_+^{\al}$ with a loss of regularity only occurring
at $s=0,\pi$ (with two-sided differentials at any order). In particular $\sigma_{\al}(\cdot)$ approximates
in the $\con^0$ sense the function $\sigma:[0,2\pi]\to\C\sim\R^2$, defined as 
$\sigma(s)=2\bar r\, e^{i(\pi-s)}$, which is a parameterization of the circle of radius 
$2\bar r$.

Taking into account Lemma~\ref{lemma-nonlocop}, the cut-locus in $V_\al$
is contained inside the set of points $Q\in \R^2$, besides $P_\eps^{\al}$, such that there exists
$(s,s',t)\in (0,\pi)\times(\pi,2\pi)\times(0,2\pi)$ for which
$Q=e^{t\tilde{X}^\al_+}\sigma_\al(s')=e^{t\tilde{X}^\al_-}\sigma_\al(s)$.

In view of applying an inverse function result for characterizing this set,
we consider the map $\Phi$ defined on
$[0,\pi]\times[\pi,2\pi]\times[0,\pi]$ by
$$
\Phi(s,s',t):=(s,e^{t\tilde{X}^{pen}_+}\sigma(s')-e^{t\tilde{X}^{pen}_-}\sigma(s)),
$$
which takes values in $\R^3$. Similarly, for $k\geq 1$ and $\al_k$ as in the proposition, we consider
the map $\Phi_k$ defined on
$[0,\pi]\times[\pi,2\pi]\times[0,\pi]$ by
$$
\Phi_k(s,s',t):=(s,e^{t\tilde{X}^{\al_k}_+}\sigma_{\al_k}(s')-e^{t\tilde{X}^{\al_k}_-}\sigma_{\al_k}(s)).
$$
Note that, since the vector fields $\tilde{X}^{\al_k}_\eps$ converge uniformly to $\tilde{X}^{pen}_\eps$ on $V_\al$, it is easy to see that $\Phi_k$ converges to $\Phi$ in the $\con^1$ norm.

For $(Pen)$, a point of the overlap curve, besides $P_\eps^{\al}$, is then identified with a triple 
$(s,s',t)\in (0,\pi)\times(\pi,2\pi)\times(0,\pi)$ such that
$\Phi(s,s',t)=(s,0,0)$. In other words, the overlap curve can be parameterized by means of  the map
$w:[0,\pi]\rightarrow \R^3$ defined implicitly by $\Phi(w(s))=(s,0,0)$, while $\ga_{pen}^o$
can be obtained as the composition of the two maps $e^{t\tilde{X}^{pen}_-}\sigma(s)$ and $w(s)$.

Similarly, we would like to define the overlap curve corresponding to $(\Sigma)_{\al_k}$,
for $k$ large enough, by means of the function $w_k$ defined by $\Phi_k(w_k(s))=(s,0,0)$.
To proceed, we will apply Theorem~\ref{seq0}. The first task consists of computing
$\det D\Phi$ along the overlap curve.

\bl\label{det0}
Along the set of triples $(s,s',t)\in (0,\pi)\times(\pi,2\pi)\times(0,\pi)$ for which
$e^{t\tilde{X}^{pen}_+}\sigma(s')=e^{t\tilde{X}^{pen}_-}\sigma(s)$, we have
$$
\det\D\Phi(s,s',t)=\frac{4\bar r(1-\bar{r}^2)\sin s'}{(1-\bar r\cos s')^2+(\bar r\sin s')^2}.
$$
\el
\noindent {\bf Proof of Lemma~\ref{det0}.} One has
$$
\det\D\Phi(s,s',t)=
\det\big((e^{t\tilde{X}^{pen}_+})_{*}\frac{d\sigma}{ds'},\tilde{X}^{pen}_+ e^{t\tilde{X}^{pen}_+}\sigma(s')
-\tilde{X}^{pen}_- e^{t\tilde{X}^{pen}_-}\sigma(s)\big).
$$
By taking into account that $\Phi(s,s',t)=0$, the previous determinant is equal to twice
the first component of $(e^{t\tilde{X}^{pen}_+})_{*}\frac{d\sigma}{ds'}$, i.e.,
$\det\D\Phi(s,s',t)=4\bar r \sin(s'-t)$. Using Eq.~\r{over-param}, one concludes.
\eproof

\vspace{0.2cm}

Observe that  $\det D\Phi\neq 0$ if $s'\neq 0,\pi$. In particular, if we consider a closed interval 
$I\subset (0,\pi)$, then the set Im$(w)_{|I}$ plays the role of the compact set $\mathcal{K}$ in 
Theorem~\ref{seq0}. All the assumptions of the theorem are then verified, and therefore
we have proved the existence of a map $w_k$ defined on $I$, satisfying
$\Phi_k(w_k(s))=(s,0,0)$ and converging uniformly to $w$.
If we define $\ga_k^o$ as the composition of the two maps $e^{t\tilde{X}^{\al_k}_-}\sigma(s)$ 
and $w(s)$, then, since $I$ was chosen arbitrarily, the proof of the theorem is complete.


\eproof

\subsubsection{Proof of Proposition~\ref{p2-th}}

With the previous notations, let $\varphi_k$ be the smooth map defined on 
$[0,\pi]\times[\pi,2\pi]\times[0,2\pi]$ by
$$
\varphi_k(s,s',t)=e^{t\tilde{X}^{\al_k}_+}\sigma_k(s')-e^{tX^{\al_k}_-}\sigma_k(s).
$$
For the rest of this paragraph, we drop the index $k$ to get lighter notations.

From the Taylor expansion of $\varphi$ around the points $(0,2\pi,0)$ and
$(\pi, \pi,0)$, we derive the asymptotic
behaviors of the cut locus close to the points $P_\eps^{\al}$, $\eps=\pm$, since that cut locus
belongs to the level set $\varphi=0$.
We will only perform computations at $(0,2\pi,0)$ since they are entirely similar at
$(\pi, \pi,0)$.

Let us call $\varphi^{(1)}$ $\varphi^{(2)}$ the two components of $\varphi$.
We use $\varphi^{(i)}_s$ to denote the partial derivative of the component $\varphi^{(i)}$
with respect to $s$ evaluated in $(0,2\pi,0)$ and we define in an analogous
way all the (multiple) partial derivatives evaluated in $(0,2\pi,0)$. Set $\tilde{s}:=s'-2\pi$.
Then, after computations, we have $\varphi^{(1)}_s=\varphi^{(1)}_{\tilde s}=\varphi^{(1)}_t=0$ and
$$
\ba{ccc}
\varphi^{(1)}_{ss}=-2\bar r+\O(\al)\,  & \varphi^{(1)}_{\tilde s\tilde s}=2\bar r+\O(\al)\,,
& \varphi^{(1)}_{tt}=2+\O(\al)\,,\\
\varphi^{(1)}_{s\tilde s}=0\,, & \varphi^{(1)}_{st}=-2\bar r+\O(\al)\,, & \varphi^{(1)}_{\tilde st}=
2\bar r+\O(\al)\,,\\
\varphi^{(2)}_s=2\bar r+\O(\al)\,, & \varphi^{(2)}_{\tilde s}=-2\bar r+\O(\al)\,, & 
\varphi^{(2)}_t=2\bar r+\O(\al)\,.
\ea
$$
We thus get
\bqn
\varphi^{(1)}(s,\tilde s,t)&=&\varphi^{(1)}_{ss}s^2+\varphi^{(1)}_{\tilde s\tilde s}
{\tilde s}^2+\varphi^{(1)}_{tt}t^2+2\varphi^{(1)}_{st}st+2\varphi^{(1)}_{\tilde st}
\tilde st+\O(|(s,\tilde s,t)|^3)\nn\\
&=&-2\bar r s^2+2r {\tilde s}^2+2t^2-4\bar r st+4r \tilde s t+\O(\al|(s,\tilde s,t)|^2)+
\O(|(s,\tilde s,t)|^3)  \,,
\eqnl{fi1}
and
\bqn
\varphi^{(2)}(s,\tilde s,t)=\varphi^{(2)}_s s+\varphi^{(2)}_{\tilde s} \tilde s+
\varphi^{(2)}_t t+\O(|(s,\tilde s,t)|^2)=2\bar r s-2r {\tilde s}-2t+
\O(\al|(s,\tilde s,t)|)+\O(|(s,\tilde s,t)|^2)\,,
\eqnl{fi2}
where, here, $\O(\cdot)$ is uniform with respect to $\al$.

Fix $\xi_0>0$ small. We are looking at the cut locus in a neighborhood of $P_\eps^{\al}$,
and thus, we can assume $|(s,\tilde s,t)|<\xi_0$ for some $\xi_0>0$.
The purpose of subsequent computations consists of expressing $\tilde{s}<0$
and $t>0$ as functions of $s$, for $0\leq s\leq \xi_0$, by using the equations
$\varphi^{(1)}=0$ and $\varphi^{(2)}=0$.


From $\varphi^{(2)}=0$, by applying the implicit function theorem, for $\xi_0$ small enough
and $|(s,\tilde s)|<\xi_0$, we get $t=h(s,\tilde s)$ with $h\in\con^1$.  Moreover, since 
$h(s,\tilde s)=\O(|(s,\tilde s)|)$ we have
\beq\label{t-ss'}
h(s,\tilde s)=\bar r s-\bar r {\tilde s}+\O(\al|(s,\tilde s)|)+\O(|(s,\tilde s)|^2).
\eeq
Consider now the map $$\phi(s,\tilde s)=\frac{\varphi^{(1)}(s,\tilde s,h(s,\tilde s))}{s-\tilde s}\,,$$
which is well defined and $\con^1$ for $s>0,\tilde s<0$. Again, it is possible to apply the 
implicit function theorem to the equation $\phi=0$, so that we get $\tilde s$ as a 
$\con^1$ function of $s$ and this gives the existence of the overlap 
curve. Moreover, by combining Eqs. \r{fi2} and \r{t-ss'}, we get the following
\beq\label{s'-s}
\phi(s,\tilde s)=s(1+\bar r)+\tilde s(1-\bar r)+\O(\al|(s,\tilde s)|)+\O(|(s,\tilde s)|^2)=0,
\eeq
and then, $|\tilde s|=\O(|s|)$.
Therefore, from this estimate and the above ones, we immediately obtain that
\beq\label{ouf1}
\tilde s=-\Big(\frac{1+\bar r}{1-\bar r}+\O(\al)\Big)s+\O(s^2)\,,\qquad
t=\Big(\frac{2\bar r}{1-\bar r}+\O(\al)\Big)s+\O(s^2)\,,
\eeq
from which we get  that the overlap curve converges uniformly to $\ga^o_{pen}$.

The proof of Proposition~\ref{p2-th} is now complete.\eproof

\brem\label{extend}
In order to prove that the overlap curve $\ga^o_k$ converges to 
$\ga^o_{pen}$ uniformly with respect to $\bar r$ in any closed interval 
$I\subset (0,1)$, it is enough to follow the lines of the proofs of Propositions~\ref{p1-th},~\ref{p2-th}   
by considering $\bar r$ as an additional variable. For instance, for the Proposition~\ref{p1-th},
one needs to define the maps $\tilde\Phi:\R^4\to\R^4$, $\tilde\Phi(\bar r,s,s',t):=
(\bar r,s,e^{t\tilde{X}^{pen}_+}\sigma(s')-e^{t\tilde{X}^{pen}_-}\sigma(s))$ (recall that 
$\sigma(\cdot)$ depends on $\bar r$) and $\tilde\Phi_k:\R^4\to\R^4$, $\tilde\Phi_k(\bar r,s,s',t):=
(\bar r,s,e^{t\tilde{X}^{\al_k}_+}\sigma_k(s')-e^{t\tilde{X}^{\al_k}_-}\sigma_k(s))$ (where 
$\al_k=\pi/(2(\bar r+k))$ and $\sigma_k(\cdot)$ depends on $\bar r$).\\
The uniformity with respect to $\bar r$ is then proved by applying Theorem~\ref{seq0} to 
$\tilde\Phi,\Tilde\Phi_k$ with $\mathcal{K}= \{ \tilde{\Phi}^{-1}(\bar r,s,0,0)\,:\ (\bar r,s)\in I\times J\}$ 
where $I\times J$ is a compact subset of $(0,1)\times(0,\pi)$.
\erem

\section{Case $r=C\al$}\label{Kal}
\subsection{Description of the minimum time front at time $k_M\pi$}
Fix $C>0$ and consider the sequence $(\al_k)$ such that $r(\al_k)=C\al_k$, $k\geq 0$.
As before, we drop the index $k$ when possible. For $\al_k$ small enough,
one deduces, from the analysis of \cite{q5}, that the south pole is not reached at time
$k_m\pi=[\frac{\pi}{2\al}]\pi$, so that the optimal front at time $k_M\pi$ is not empty. 
The next result provides a description of the extremal front
at time $k_M\pi$.

\bl\label{lem-est2}
Define the planar curve ${\cal L}:[0,2\pi]\rightarrow \R^2$ by
\beq\label{LL0}
{\cal L}(s)=\bpm \cs (-2C+\pi \sis^2/2) \\ \sis(\pi+2C-\pi\sis^2/2)\epm.
\eeq
Then, for $s\in [0,\pi]$, we have
\beq\label{f+02}
\mathcal{E}^+(\al,s)=(\al^2{\cal L}(s),-1)^T +\O(\al^3), \eeq
and \beq\label{df+02} \derp{s}\mathcal{E}^+(\al,s)=(\al^2\dert{s}{\cal L}(s),0)^T
+\O(\al^3).
\eeq
At $s=0$ and $s=\pi$, the derivatives
are only one-sided, i.e. as $s>0$ tends to zero and $s<\pi$ tends to
$\pi$.

Similarly, we have, for $s\in [0,\pi]$,
\beq\label{f-02}
\mathcal{E}^-(\al,k_M\pi,s)=(\al^2{\cal L}(s+\pi),-1)^T +\O(\al^3), \eeq
and \beq\label{df-02} \derp{s}\mathcal{E}^-(\al,k_M\pi,s)=(\al^2\dert{s}{\cal L}(s+\pi),0)^T
+\O(\al^3),
\eeq
with one-sided derivatives at $s=0$ and $s=\pi$.
\el

\noindent {\bf Proof of Lemma~\ref{lem-est2}.} This is immediate from Proposition~\ref{pop0} applied in the case
$r(\al)=C\al$.
\eproof

For $C<\pi/4$, consider $\theta_d\in (0,\pi/2)$ with $\sin(\theta_d)=2\sqrt{C/\pi}$.
The curve ${\cal L}(s)$ has two double points
$D^+={\cal L}(s_1^+)={\cal L}(s_2^+)$, with $s_1^+=\theta_d$ and $s_2^+=\pi-\theta_d$,
and $D^-={\cal L}(s_1^-)={\cal L}(s_2^-)$, with $s_1^-=\pi+\theta_d$ and
$s_2^-=2\pi-\theta_d$. It also has four cuspidal points $Cp^\eps_i$, $i=1,2$ and
$\eps=\pm$, corresponding to the values $s=s_{cusp,i}^\eps$, where $\sin^2 s=\frac{2+4C/\pi}3$.

Finally, let $\sigma$ be the closed Jordan curve defined as
the restriction of ${\cal L}(s)$ to $[0,s_1^+]\cup[s_2^+,s_1^-]\cup [s_2^-,2\pi]$.
If $C>\pi/4$, we simply define $\sigma$ to be ${\cal {L}}$.

At the light of the previous result, we get that $\F(\al,k_M\pi)$, the extremal front at time $k_M\pi$, is contained
inside a neighborhood $W_\al$ of the south pole of order $\O(\al^2)$ neighborhood of the south pole.  
Therefore, in order to
understand the shape of the optimal synthesis inside $W_\al$, we must rescale the whole
problem by $N_\al$, the linear
mapping from to $\R^3$ onto $\R^2$ defined as the composition of the orthogonal projection
$(x_1,x_2,x_3)\mapsto (x_1,x_2)$ followed by the dilation by $1/\al^2$.

For $x\in W_\al$, we first consider $(\Lambda)_{\al}$, the image of $(\Sigma)$ by $N_\al$,
i.e. $(\Lambda)_{\al}$
is the planar control system given by
\bqn (\Lambda)_{\al}:\quad\left\{ \ba{l}
  \dot{z_1}=-\ca  z_2, \\
  \dot{z_2}=\ca z_1 +u\frac{\sa}{\al^2} \sqrt{1-\al^4\|z\|^2}. \ea \right.
\eqnl{f-lam}
Let ${\cal L}_\al$ be the image of $\F(\al,k_M\pi)$ by $N_\al$.
From Lemma~\ref{lem-est2}, ${\cal L}_\al$ converges to ${\cal L}$ in the
$C^1$ topology. It is clear that, for $C>\pi/4$,
${\cal L}_\al:[0,2\pi]\rightarrow \R^2$ is homeomorphic to $e^{is}$, $s\in [0,2\pi]$.
In the case where $C<\pi/4$, the next lemma shows that, for $\al$
small enough, ${\cal L}_\al$ has the same shape as ${\cal L}$.

\ppotR{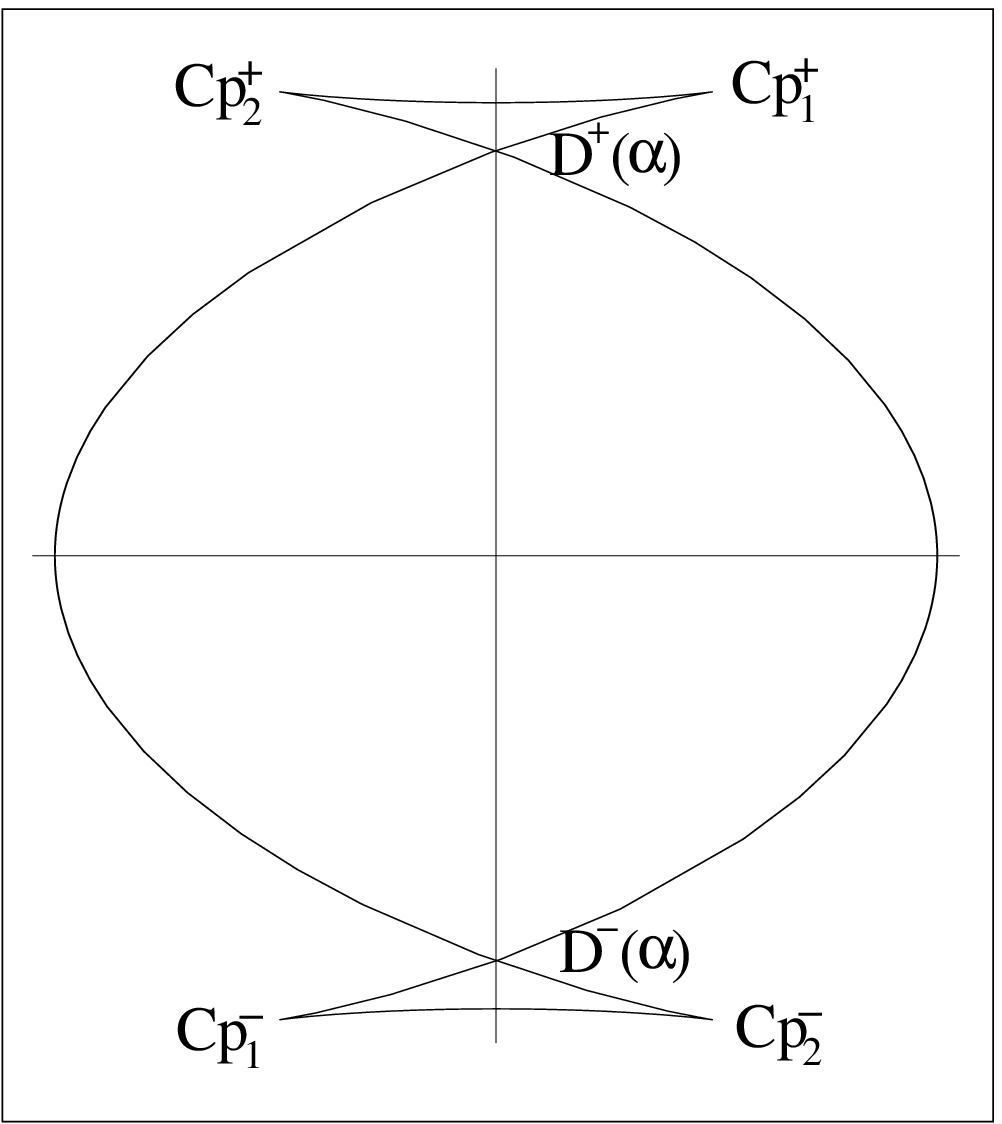}{Graph of the function ${\cal L}_\al$ for $C<\pi/4$}{7}

\bl\label{shape}
If $C<\pi/4$, then ${\cal L}_\al$ is described by the following picture,
where $Cp^\eps_i(\al)={\cal L}_\al(s_{cusp,i}^\eps(\al))$, $i=1,2$ and $\eps=\pm$,
are cuspidal points and $D^\eps(\al)$
are double points with
\beq\label{ddp}
D^+(\al)={\cal L}_\al(s_1^+(\al))={\cal L}_\al(s_2^+(\al)),\qquad
D^-(\al)={\cal L}_\al(s_1^-(\al))={\cal L}_\al(s_2^-(\al)),
\eeq
where $s_{cusp,i}^\eps(\al)$ and $s_i^\eps(\al)$ tend respectively to
$s_{cusp,i}^\eps$ and $s_i^\eps$ as $\al$ tends to zero, for
$i=1,2$ and $\eps=\pm$. For $\al$ small enough, set $\sigma_\al$, the
closed curve defined as the restriction of ${\cal L}_\al(s)$ to
$[0,s_{1}^+(\al)]\cup[s_{2}^+(\al),s_{1}^-(\al)]\cup [s_{2}^-(\al),2\pi]$.
Then, it is a Jordan curve.
\el

\noindent {\bf Proof of Lemma~\ref{shape}.} For $i=1,2$ and $\eps=\pm$, the existence of the
cuspidal points $Cp^\eps_i(\al)$ is obtained by applying the implicit function
theorem to the equation $DL(s,\al)=0$, where the function $DL(s,\al)$ is defined by
$$
DL(s,\al):=\dert{s}{\cal {L}}_{\al}(s),
$$
in the neighborhood of each $(s_{cusp,i}^\eps,0)$. We have
$$
\partial_sDL(s_{cusp,i}^\eps,0)=\frac{d^2}{ds^2}{\cal {L}}(s_{cusp,i}^\eps)\neq 0
$$
and we conclude. The uniqueness of these four points,
on $[0,2\pi]$, is trivial since $DL(s,\al)=\dert{s}{\cal {L}}(s)+\O(\al)$.

Similarly, for $\eps=\pm$, the existence of the double points $D^\eps(\al)$
follows after applying the implicit function theorem to the equation $DP(s,s',\al)=0$, 
where the function $DP(s,s',\al)$ is defined by
$$
DP(s,s',\al)={\cal {L}}_{\al}(s)-{\cal {L}}_{\al}(s'),
$$
in the neighborhood of each $(s_1^\eps,s_2^\eps,0)$. For the uniqueness, we
proceed as before.
\eproof

In the case $C>\pi/4$, we also define $\sigma_\al$ to be equal to ${\cal L}_\al$.
As a consequence, we are able to characterize $OF(\al,k_M\pi)$, the minimum time front
at time $k_M\pi$ when $C\neq \pi/4$.

\bp\label{OF1}
For $\al$ small enough and $C\neq \pi/4$, the minimum time front at time $k_M\pi$,
$OF(\al_k,k_M\pi)$ is equal to $\tilde{\sigma}_\al$, the inverse image on $S^2$,
by $N_\al$, of $\sigma_\al$.
\ep

\brem
As a consequence, we deduce that, for $C>\pi/4$ and $\al$ small enough, the optimal synthesis
between $\F(\al,(k_M-1)\pi)$ and $\F(\al,k_M\pi)$ is simply given by
the extremal flow whereas, for $C<\pi/4$, there is a loss of
optimality along certain extremal curves starting at
$\F(\al,(k_M-1)\pi)$ before reaching $\F(\al,k_M\pi)$. The values of
$s$ corresponding to such curves can be deduced from the previous
characterizations of $\F(\al,k_M\pi)$ and $OF(\al_k,k_M\pi)$.
\erem

\noindent {\bf Proof of Proposition~\ref{OF1}.} Recall that $OF(\al_k,k_M\pi)$ is a
piecewise $\con^1$ submanifold of $\F(\al,k_M\pi)$. As in the proof of
Proposition~\ref{pop2}, the result to establish is a consequence of the fact that 
$\sigma_\al={\cal  L}_\al:[0,2\pi]\rightarrow \R^2$ is homeomorphic to $e^{is}$, $s\in
[0,2\pi]$ and it can be achieved by means of simple topological arguments.

In the case $C<\pi/4$, $\sigma_\al$ is a piecewise $\con^1$ Jordan curve
homeomorphic to $e^{is}$, $s\in [0,2\pi]$. A simple topological argument yields
the conclusion.\eproof

\subsection{Limit of the synthesis}

It remains to describe the limiting dynamics close to the south pole.
In order to take the limit, as $\al$ tends to zero, in
$(\Lambda)_\al$, one must reparameterize by the time $\al t$. The
limit is then given by the control system
$$
(\Lambda):\qquad \left\{ \ba{l}
  \dot{z_1}=0, \\
  \dot{z_2}=u. \ea \right.
$$
We now describe the optimal synthesis for the limit problem, i.e. for
the problem of reaching in minimum time every point inside $\sigma$
along $(\Lambda)$ and starting from $\sigma$. Because of the
symmetries of $\sigma$ and because the tangent vector to $\sigma$ is
vertical only at $s=0$ and $s=\pi$, there exists a unique overlap
curve $(Seg)_C$, defined as the segment of the $z_1$-axis between the points
$(-2C,0)$ and $(2C,0)$. Above it, the input $u$ takes the constant
value $-1$ and, below that overlap curve, the constant value $1$.
Integral curves are clearly vertical lines.

We next intend to prove that the optimal synthesis consisting
of reaching in minimum time every point inside $\sigma_\al$
along $(\Lambda)_\al$ and starting from $\sigma_\al$ converges
to the previous synthesis in the following sense.

\bt\label{the3}
Assume that $C\neq \pi/4$. As $\al$ tends to zero, the time
optimal synthesis associated to $(\Lambda)_{\al}$ inside
$\sigma_\al$ is characterized by an overlap curve $(Seg)_C^\al$,
converging to $(Seg)_C$ in the $C^0$ topology, and, above
$(Seg)_C^\al$, the control $u$ takes the constant value $-1$ and
below $(Seg)_C^\al$, it is equal to $1$. Moreover, there exist only two
time optimal trajectories reaching the origin and, in the case $C<\pi/4$,
these trajectories start from $D^\eps_\al$, $\eps=\pm$, the
double points of ${\cal {L}}_\al$.
\et

\noindent {\bf Proof of Theorem~\ref{the3}.}
Fix $C\neq \pi/4$. We first notice that, for $\al$ small enough, there are not
switching curves inside $\sigma_\al$. Therefore, the cut-locus may only occur
as images by $N_\al$ of points $M\in S^2$ such that
$M=e^{\frac{t}{\al}X_-}\tilde{\sigma}(s)=e^{\frac{t}{\al}X_+}\tilde{\sigma}(s')$
for $t\in [0,\frac{2\pi}{\al}]$, $s\in [0,\pi]$ and $s'\in [\pi,2\pi]$.
Proceeding exactly as in the proof of Theorem~\ref{the1}, we apply inverse function
arguments first in neighborhoods of $\sigma_\al(0)$ and $\sigma_\al(\pi)$ and second
in a region enclosed by $\sigma_ \al$ excluding such neighborhoods. It is then easy to
determine  the values of the input $u$ in each connected component of the region
enclosed by $\sigma_ \al$ minus $(Seg)_C^\al$.

By a continuity argument, it is clear that there exist only two time
optimal trajectories reaching the origin, one above $(Seg)_C^\al$ and
one below. Finally, suppose that $C<\pi/4$. In that case, it was
proved in \cite{q5} that the only extremals starting at a point ${\cal
  {L}}_\al(s)$ and reaching the origin from above the overlap curve
$(Seg)_C^\al$ correspond to values of $s$ verifying one of the following three
possibilities as $\al$ tends to zero: (a) $s$ tends
to zero, (b) $s$ tends to $\pi/2$, (c) ${\cal {L}}_\al(s)$
is a double point also associated to $s'=v(s)-s$. In view of what precedes, only
possibility (c) is allowed for optimality. Theorem~\ref{the3} is proved.\eproof

\brem
As a consequence of the previous argument and from the results of \cite{q5}, 
we get that, for $\al$ small enough and $C<\pi/4$,
$$
s_2^+(\al)=v(s_1^+(\al))-s_1^+(\al),\
s_2^-(\al)=2\pi+ v(s_1^-(\al)-\pi)-s_1^-(\al),
$$
where $s_i^\eps(\al)$, $i=1,2$ $\eps=\pm$, were defined in \r{ddp}.
\erem

\section{Case $r(\al)=0$}\label{r=0}
We assume here that $r(\al)=0$, i.e. $\al_k=\frac{\pi}{2k}$ for $k\geq 1$.
From  Proposition~\ref{pop0},
we know that the extremal front at time
$([\frac{\pi}{2\al}]-1)\pi=\frac{\pi}{2\al}-\pi$, encloses the south pole,
is optimal and is approximately (in the $\con^1$ sense) a circle of radius $2r(\al)\al$
around the south pole. Moreover, at time $[\frac{\pi}{2\al}]\pi$, we know that
the extremal front must contain the south pole and is equal, up to $\O(\al^3)$,
to $(\al^2{\cal{L}},-1)^T$ given in \r{f+02} and \r{f-02} with $C=0$.
In that case, the minimum time front reduces to the south pole.

In this case it is interesting to consider the synthesis starting from the extremal front at time 
$(k_M-1)\pi$ and it is natural to compare it with the synthesis of the linear pendulum
studied in Section~\ref{s-rwlp} and corresponding to
$\rho=2$. Let us first describe briefly that synthesis. Let $D_2$ and $C_2$
be the disc and the circle centered at the origin and of radius $2$ respectively.
The overlap curve inside $D_2$ coincides with the switching
curves and with the trajectories, corresponding to $u=\pm 1$,
connecting the points $(\pm 2,0)$ to the origin. In particular, it means that
an optimal trajectory of the synthesis starting at any point $P\in C_2$
reaches the origin, and thus, there exists an infinite number of optimal trajectories
from $C_2$ to the origin.

For $\al>0$ and $r(\al)=0$, the situation is rather different.
Let us first define $\tilde{F}(\al,(k_M-1)\pi)$ to be the image
of $\F(\al,(k_M-1)\pi)$ by $M_\al$. Then, for $\al$ small enough,
it was shown in \cite{q5}, that the only optimal trajectories starting from
$\tilde{F}(\al,(k_M-1)\pi)$ and reaching the origin are those starting at
$P^\al_+$ and $P^\al_-$. Let us refer to them as $\ga^+$ and $\ga^-$.
Therefore, in the case  $r(\al)=0$, the synthesis
for $\al>0$ is rather different than the synthesis of the limit candidate
when $\al$ tends to zero. It is a clear indication that the case $r(\al)=0$
is more delicate than the cases $r(\al)$ positive constant or $r(\al)=C\al$.
However we are still able to give a partial description of the limit synthesis
as  the next proposition shows. 

\bp\label{Pr=0}
Assume that $r(\al)=0$ and $\al$ is small enough. Then the switching curve
$C^+_{k_M}$ (resp. $C^-_{k_M}$) is optimal for some interval $[0,s(\al)]$,
$s(\al)<\pi$, and it is above (resp. below) $\ga^+$ (resp. $\ga^-$)
as long as it is optimal. Moreover, we have
\beq\label{lolo}
\lim_{ \al\rightarrow 0, r(\al)=0}s(\al)=\bar{s}:=\arccos{\sqrt{1/3}}.
\eeq
\ep
\noindent {\bf Proof of Proposition~\ref{Pr=0}.}
We only provide an argument for $C^+_{k_M}$, being the other case analogous. 
To prove the first statement of the proposition, we reason by contradiction. If the switching 
curve is not optimal on any interval $[0,\tau]$, $\tau>0$, we get
the existence of an optimal trajectory starting at
$\F(\al,(k_M-1)\pi)$ above $P^+_\al$ and reaching the origin, which is
equal to the concatenation of an integral curve of $X_-$ and a piece
of $\ga^+$.  Therefore, an optimal integral curve of $X_-$, starting
above $\ga^+$, must either switch or lose optimality before reaching
$\ga^+$. If the second possibility occurs, we must have an overlap
i.e., at that point an optimal integral curve of $X_+$ arrives. 
Close to $P^+_\al$, the latter would imply that the
optimal integral curve of $X_+$ starts at $\F(\al,(k_M-1)\pi)$ above
$P^+_\al$. This is impossible because, from every point of
$\mathcal{E}^+(\al,(k_M-1)\pi)$, the value of the optimal control is $-1$.
Let $s(\al)\leq \pi$ be the first value of $s$ for which $C^+_{k_M}$ ceases to
be optimal. Define
$$
H(s):=\det\big(X_+(C^+_{k_M}(s)), \frac{dC^+_{k_M}}{ds}(s),C^+_{k_M}(s)\big),
$$
for $s\in [0,\pi]$. Then, $s(\al)$ is the smallest solution in $(0,\pi]$
of $H(s)=0$. It is easy
to see that $H$ must take the value zero before $\pi$. We deduce that
$s(\al)<\pi$. By taking the asymptotic expansion of the previous
expression as $\al$ tends to zero, we get
$$
H(s)=\frac{\pi}4\sis \al^3(1+3\cos(2s)+\O(\al)).
$$
Then $s(\al)$ must converge to $\bar{s}$ as $\al$ tends to zero, the smallest
solution in $[0,\pi]$ of $1+3\cos(2s)=0$.
\eproof


\ppotR{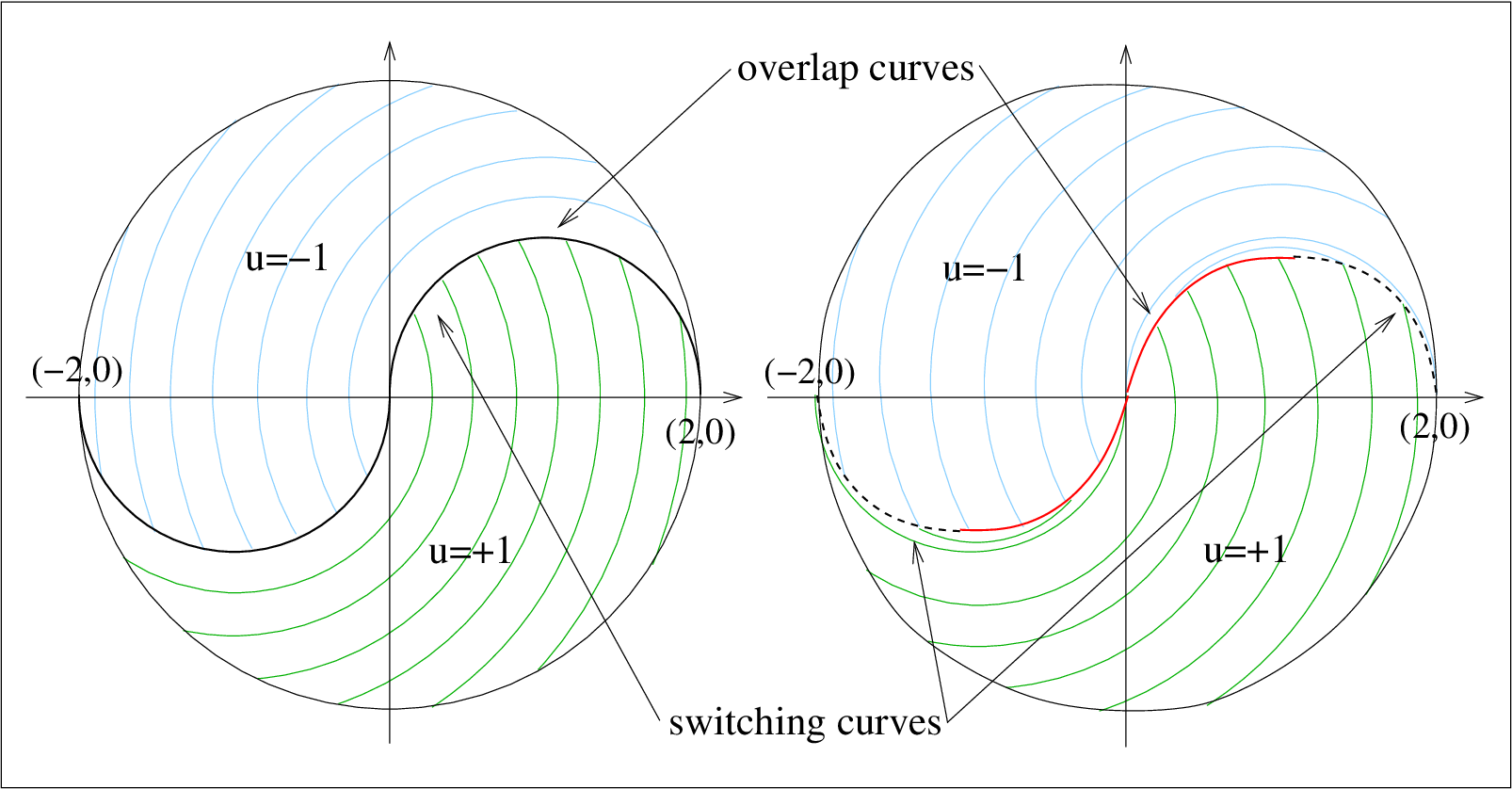}{Comparison between the optimal synthesis for the linear pendulum
and the optimal synthesis on the bottom of the sphere in the case $\ra=0$}{12}


\section{Appendix}
The following version of the inverse function theorem is used in the argument
of Proposition~\ref{p1-th}.

\bt\label{seq0}
Let $\Phi:\R^n\to\R^n$ be a $\con^1$ map and ${\cal {K}}\subset \R^n$ a compact set such
that $\Phi_{|{\cal {K}}}:{\cal {K}}\to\Phi({\cal {K}})$
is bijective and the differential $D\Phi(x)$ is invertible for $x\in K$.
Then, there exists an open neighborhood
$U\supset {\cal {K}}$ such that $\Phi_{|U}$ is
a $\con^1$ diffeomorphism.

Let now $(\Phi_k)_{k\geq 1}$ be a sequence of $\con^1$ maps converging
in the $\con^1_{loc}$ sense to $\Phi$. Then, for every open set
$\tilde{U}$ with closure included in $U$, there exists $\bar k$ such
that, for every $k\geq\bar k$, $\Phi_{k_{|\tilde{U}}}$ is a $\con^1$
diffeomorphism and, for every compact subset $\tilde{{\cal {K}}}$ of
$\tilde{U}$, $\Phi(\tilde{{\cal {K}}})\subset \Phi_k(\tilde{U})$ and
$\displaystyle{\lim_{k\to\infty} \Phi_k^{-1}(v)=\Phi^{-1}(v)}$
uniformly with respect to $v\in \Phi(\tilde{{\cal {K}}})$.
\et

\proof
Let us define, for $k\geq 0$, the following open neighborhoods of ${\cal {K}}$
$$A:=\{x\in \R^n\,:\,\det D\Phi(x)\neq 0\}
\qquad A_k:=\cup_{x\in {\cal {K}}}B\Big(x,\frac{1}{k}\Big)\cap
A\,.$$ In view of the inverse function theorem, in order to conclude
the proof of the first part, it is enough to show that for $k$ large
enough the restriction $\Phi_{|A_k}$ is one-to-one.

We argue by contradiction. Let $x_k\neq y_k\in A_k$ such that
$\Phi(x_k)=\Phi(y_k)\quad \forall k$. Then, up to extractions of subsequences,
we can assume that the two sequences converge to $\bar x$ and  $\bar y$
respectively. Since  $\bar x, \bar y\in\cap_k A_k={\cal {K}}$ and
$\Phi(\bar x)=\Phi(\bar y)$, we deduce that $\bar x=\bar y$. However,
since $\det D\Phi(\bar x)\neq 0$, we have that $\Phi$ is bijective in a neighborhood of
$\bar x$, which contradicts the assumption $\Phi(x_k)=\Phi(y_k)$ for $k$
large enough.

The proof of the second part is similar. First, fix a subset
$\tilde{U}$ of $U$. By the uniform convergence of $D\phi_k$ to $D\phi$
on every compact subset of $U$, we get $\det D\Phi_k(x)\neq 0$ for
every $x\in\tilde{U}$ and $k$ large enough. We also obtain that
$\Phi_k$ is one-to-one with the same argument as above. For the remaining results
to establish, they simply follow from the uniform convergence of $\Phi_k$
to $\Phi$ on every compact subset of $U$.\eproof

\end{document}